\documentclass[11pt]{amsart}

\usepackage{eufrak}
\usepackage{euscript}
\usepackage{amsmath}
\usepackage{amsfonts}
\usepackage{amssymb}

\usepackage{amsxtra}
\usepackage{amscd}
\usepackage{amsthm}

\usepackage{graphicx}
\usepackage{latexsym}

\newtheorem{Definition}{Definition}[section]
\newtheorem{Proposition}{Proposition}[section]
\newtheorem{Lemma}{Lemma}[section]
\newtheorem{Theorem}{Theorem}[section]
\newtheorem{Corollary}{Corollary}[section]
\newtheorem{Remark}{Remark}[section]
\newtheorem{Example}{Example}[section]

\newcommand{\nc}{\newcommand}
\newcommand{\on}{\operatorname}

\nc{\BA}{{\mathbb{A}}}
\nc{\BC}{{\mathbb{C}}}
\nc{\BE}{{\mathbb{E}}}
\nc{\BK}{{\mathbb{K}}}
\nc{\BM}{{\mathbb{M}}}
\nc{\BN}{{\mathbb{N}}}
\nc{\BP}{{\mathbb{P}}}
\nc{\BR}{{\mathbb{R}}}
\nc{\BZ}{{\mathbb{Z}}}
\nc{\BS}{{\mathbb{S}}}

\nc{\CA}{{\mathcal{A}}}
\nc{\CB}{{\mathcal{B}}}
\nc{\CC}{{\mathcal{C}}}
\nc{\CK}{{\mathcal{K}}}
\nc{\CE}{{\mathcal{E}}}
\nc{\CF}{{\mathcal{F}}}
\nc{\CG}{{\mathcal{G}}}
\nc{\CI}{{\mathcal{I}}}
\nc{\CJ}{{\mathcal{J}}}
\nc{\CL}{{\mathcal{L}}}
\nc{\CM}{{\mathcal{M}}}
\nc{\CN}{{\mathcal{N}}}
\nc{\CO}{{\mathcal{O}}}
\nc{\CP}{{\mathcal{P}}}
\nc{\CQ}{{\mathcal{Q}}}
\nc{\CR}{{\mathcal{R}}}
\nc{\CS}{{\mathcal{S}}}
\nc{\CT}{{\mathcal{T}}}
\nc{\CU}{{\mathcal{U}}}
\nc{\CV}{{\mathcal{V}}}
\nc{\CW}{{\mathcal{W}}}
\nc{\CX}{{\mathcal{X}}}
\nc{\CY}{{\mathcal{Y}}}
\nc{\CZ}{{\mathcal{Z}}}

\nc{\fa}{{\mathfrak{a}}}
\nc{\fb}{{\mathfrak{b}}}
\nc{\fg}{{\mathfrak{g}}}
\nc{\fh}{{\mathfrak{h}}}
\nc{\fj}{{\mathfrak{j}}}
\nc{\fk}{{\mathfrak{k}}}
\nc{\fl}{{\mathfrak{l}}}
\nc{\fm}{{\mathfrak{m}}}
\nc{\fn}{{\mathfrak{n}}}
\nc{\fo}{{\mathfrak{o}}}
\nc{\fu}{{\mathfrak{u}}}
\nc{\fp}{{\mathfrak{p}}}
\nc{\fr}{{\mathfrak{r}}}
\nc{\fs}{{\mathfrak{s}}}
\nc{\ft}{{\mathfrak{t}}}

\nc{\fA}{{\mathfrak{A}}}
\nc{\fB}{{\mathfrak{B}}}
\nc{\fD}{{\mathfrak{D}}}
\nc{\fE}{{\mathfrak{E}}}
\nc{\fF}{{\mathfrak{F}}}
\nc{\fG}{{\mathfrak{G}}}
\nc{\fK}{{\mathfrak{K}}}
\nc{\fL}{{\mathfrak{L}}}
\nc{\fM}{{\mathfrak{M}}}
\nc{\fN}{{\mathfrak{N}}}
\nc{\fP}{{\mathfrak{P}}}
\nc{\fU}{{\mathfrak{U}}}
\nc{\fV}{{\mathfrak{V}}}
\nc{\fZ}{{\mathfrak{Z}}}

\nc{\bb}{{\mathbf{b}}}
\nc{\bc}{{\mathbf{c}}}
\nc{\bd}{{\mathbf{d}}}
\nc{\be}{{\mathbf{e}}}
\nc{\bj}{{\mathbf{j}}}
\nc{\bn}{{\mathbf{n}}}
\nc{\bp}{{\mathbf{p}}}
\nc{\bq}{{\mathbf{q}}}
\nc{\bu}{{\mathbf{u}}}
\nc{\bv}{{\mathbf{v}}}
\nc{\bx}{{\mathbf{x}}}
\nc{\bs}{{\mathbf{s}}}
\nc{\by}{{\mathbf{y}}}
\nc{\bw}{{\mathbf{w}}}
\nc{\bA}{{\mathbf{A}}}
\nc{\bK}{{\mathbf{K}}}
\nc{\bB}{{\mathbf{B}}}
\nc{\bC}{{\mathbf{C}}}
\nc{\bD}{{\mathbf{D}}}
\nc{\bH}{{\mathbf{H}}}
\nc{\bM}{{\mathbf{M}}}
\nc{\bN}{{\mathbf{N}}}
\nc{\bP}{{\mathbf{P}}}
\nc{\bV}{{\mathbf{V}}}
\nc{\bW}{{\mathbf{W}}}
\nc{\bX}{{\mathbf{X}}}
\nc{\bZ}{{\mathbf{Z}}}
\nc{\bS}{{\mathbf{S}}}

\nc{\sA}{{\mathsf{A}}}
\nc{\sB}{{\mathsf{B}}}
\nc{\sC}{{\mathsf{C}}}
\nc{\sD}{{\mathsf{D}}}
\nc{\sF}{{\mathsf{F}}}
\nc{\sK}{{\mathsf{K}}}
\nc{\sM}{{\mathsf{M}}}
\nc{\sO}{{\mathsf{O}}}
\nc{\sS}{{\mathsf{S}}}
\nc{\sZ}{{\mathsf{Z}}}

\nc{\tZ}{{\tilde{\mathbb Z}}}

\nc{\bigimlat}{{\Sigma(\cowt)}}
\nc{\imlat}{{\sigma(\cowt)}}
\nc{\imgp}{{\EuScript L}}

\nc{\Aut}{\on{Aut}}
\nc{\Inn}{\on{Inn}}
\nc{\red}{{\on{red}}}
\nc{\loc}{{\on{loc}}}
\nc{\pos}{{\on{pos}}}

\nc{\ch}{\check}

\nc{\R}{\Bbb R}
\nc{\C}{\Bbb C}

\newcommand{\Spec}{\on{Spec}}

\nc{\Hom}{\on{Hom}}
\nc{\Comod}{\on{Comod}}
\nc{\GL}{\on{GL}}
\nc{\Gr}{\on{Gr}}
\newcommand{\affgr}{\on{Gr}}
\newcommand{\raffgr}{\on{Gr}_{\Bbb R}}

\nc{\eH}{{\EuScript H}}
\nc{\eG}{{\EuScript G}}

\nc{\chG}{{\check G}}
\nc{\chT}{{\check T}}

\nc{\sP}{{\mathsf{P}}}
\nc{\sQ}{{\mathsf{Q}}}
\nc{\Vect}{{\on{\mathsf Vect}}}

\nc{\risom}{\stackrel{\sim}{\rightarrow}}
\nc{\lisom}{\stackrel{\sim}{\leftarrow}}

\nc{\mv}{S^\nu}
\nc{\rmv}{{S}^\nu_\R}

\nc{\nadgr}{\mathsf{adGr}^{(n)}}
\nc{\adgr}{\mathsf{adGr}}

\nc{\twt}{{\on F}^\nu}
\nc{\swt}{{\on F}^\nu_\R}

\nc{\Lie}{\on{Lie}}

\newcommand{\jet}{{\mathcal J}}

\newcommand{\catd}{{\mathbf D}}

\newcommand{\tors}{\mathcal F}
\newcommand{\sh}{\mathcal F}

\newcommand{\sraffgr}{{\Gr}^{(\sigma)}_\R}

\newcommand{\alg}{A}

\newcommand{\morph}{\mbox{Hom}}
\newcommand{\spec}{\mbox{Spec}}

\newcommand{\co}{\mathcal{O}}
\newcommand{\ck}{\mathcal{K}}

\newcommand{\cor}{\mathcal{O}_\R}
\newcommand{\ckr}{\mathcal{K}_\R}

\newcommand{\cowt}{{\Lambda_T}}

\newcommand{\dcowt}{{\Lambda^+_T}}

\newcommand{\rcowt}{{\Lambda_S}}
\newcommand{\rdcowt}{{\Lambda^+_S}}

\newcommand{\blg}{\Omega G_c}

\newcommand{\bls}{\Omega X_c}

\newcommand{\B}{\mathcal{B}}
\newcommand{\orb}{\mathcal O}

\newcommand{\deltat}{\delta^\tau}

\newcommand{\tepsilon}{\tilde\epsilon}
\newcommand{\teta}{\tilde\eta}
\newcommand{\ttheta}{\tilde\theta}
\newcommand{\thetat}{\theta^\tau}
\newcommand{\tthetat}{\tilde\theta^\tau}

\newcommand{\gcore}{C^\lambda}

\newcommand{\korb}{{\mathcal O}_K^\lambda}
\newcommand{\rorb}{{\mathcal O}_{\R}^\lambda}

\newcommand{\lieg}{\EuFrak g}

\newcommand{\liegr}{{\EuFrak g}_\R}

\newcommand{\pC}{\Bbb C^\times}

\newcommand{\prim}{\dot{\leq}}

\newcommand{\tphi}{\tilde\phi}

\newcommand{\proj}{\Bbb P}
\newcommand{\z}{\Bbb Z}

\newcommand{\lockc}{K^{{loc}}_c}

\newcommand{\liek}{\EuFrak k}

\newcommand{\Rp}{\Bbb R_{>0}}
\newcommand{\Rnn}{\Bbb R_{\geq 0}}

\newcommand{\taffgr}{\Gr^{(2)}}

\newcommand{\cA}{\mathcal{A}}
\newcommand{\poly}{\C[t]}
\newcommand{\lpoly}{\C[t,t^{-1}]}
\newcommand{\coc}{\mathcal{O}}
\newcommand{\ckc}{\mathcal{K}}

\newcommand{\ipoly}{\C[t^{-1}]}
\newcommand{\ipolyr}{\R[t^{-1}]}

\newcommand{\gstrat}{{S^\lambda}}
\newcommand{\gdstrat}{{T^\lambda}}

\newcommand{\xstrat}{{P^\lambda}}
\newcommand{\xdstrat}{{Q^\lambda}}

\newcommand{\xcore}{{B^\lambda}}

\newcommand{\rgcore}{C^\lambda_\R}
\newcommand{\rgstrat}{S_\R^\lambda}
\newcommand{\rgdstrat}{T^\lambda_\R}

\newcommand{\srgdstrat}{\tilde T^\lambda_\R}

\newcommand{\thaffgr}{\on{GR}}

\newcommand{\corb}{\CO^\lambda_c}

\newcommand{\ad}{{\on{ad}}}

\begin{document}



\title[Matsuki correspondence]{Matsuki correspondence\\ for the affine Grassmannian}

\author{\sc David Nadler }
\address{Department of Mathematics\\ 
University of Chicago\\ Chicago, IL 60637}
\email{nadler@math.uchicago.edu}

\begin{abstract}
Let $\CK=\BC((t))$ be the field of formal Laurent series,
and let $\CO=\BC[[t]]$ be the ring of formal power series.
In this paper, we present a version of the Matsuki correspondence
for the affine Grassmannian $\affgr=G(\CK)/G(\CO)$ 
of a connected reductive complex algebraic
group $G$. Our main statement is an anti-isomorphism between the orbit 
posets of two subgroups of $G(\CK)$ acting on $\affgr$.
The first subgroup is the polynomial loop group $LG_\R$ of a real form $G_\R$
of $G$; the second is the loop group $K(\CK)$
of the complexification $K$
of a maximal compact subgroup $K_c$ of $G_\R$.
The orbit poset itself turns out to be simple to describe.
\end{abstract}

\maketitle


\section{Introduction}\label{smatintro}

\subsection{Matsuki correspondence}\label{ssmat}
We begin by recalling the form of the Matsuki correspondence
that is relevant to what follows.
Let $G$ be a connected reductive complex algebraic group. 
Let $G_\R$ be a real form of $G$,
and let $\theta$ denote the conjugation of $G$ with respect to $G_\R$.
Choose a Cartan involution $\delta$ of $G$ which commutes with $\theta$, 
and let $\eta$ denote the composition of $\theta$ and $\delta$.
The three commuting involutions $\theta,\delta$, and $\eta$ provide a lattice of subgroups
$$ \begin{matrix}  
 & & G & & \\ 
 & \nearrow & \uparrow & \nwarrow & \\
 G_\R && G_c && K  \\
 & \nwarrow & \uparrow & \nearrow & \\
 & & K_c & & 
\end{matrix}
$$
where $G_c$ is the fixed points of $\delta$, $K$ is the fixed points of $\eta$,
and $K_c$ is the fixed points of $\theta, \delta$, and $\eta$ simultaneously.
As $\delta$ is a Cartan involution, $G_c$ is a maximal compact subgroup of $G$,
and $K_c$ is a maximal compact subgroup of both $G_\R$ and $K$.

Let $\B$ denote the flag variety of $G$.
The groups $G_\R$ and $K$ act on $\B$ with finitely many orbits.
In general, the set of orbits of a group acting on a space
comes with a natural partial ordering: orbits $\orb_1$, $\orb_2$
satisfy
${\mathcal O_1}\leq {\mathcal O_2}$ if and only if 
${\mathcal O_1}$ is in the closure of ${{\mathcal O}}_2$.
The Matsuki correspondence is an order-reversing isomorphism 
between the orbit posets of $G_\R$ and $K$ acting on $\B$.

\begin{Theorem}[\cite{M}]\label{tmat}
For each $G_\R$-orbit $\orb_\R$ in $\B$, there exists a unique $K$-orbit $\orb_K$
such that $\orb_\R\cap\orb_K$ is a single $K_c$-orbit.
This correspondence provides an order-reversing isomorphism
from the poset $G_\R\backslash \B$ to the poset $K\backslash \B$.

\end{Theorem}

\begin{Example}\label{exmat}
The following lattice of groups is of the type described above
$$ \begin{matrix} 
 & & \on{SL}_2(\C) & & \\ 
 & \nearrow & \uparrow & \nwarrow & \\
 \on{SL}_2(\R) && \on{SU}_2 && \on{SO}_2(\C)  \\
 & \nwarrow & \uparrow & \nearrow & \\
 & & \on{SO}_2(\R) & & 
\end{matrix}
$$

Under the standard identification of the flag variety of $\on{SL}_2(\C)$ 
with the complex projective line $\C\Bbb P^1=\C\cup \{ \infty \}$,
the $\on{SL}_2(\R)$-orbits are the upper and lower halfplanes,
and the real projective line $\R\Bbb P^1=\R\cup\{\infty\}$.
The dual $\on{SO}_2(\C)$-orbits are the points $\{i\}$ and $\{-i\}$,
and their complement $\C\Bbb P^1\setminus \{i, -i\}$.
The $\on{SO}_2(\R)$-orbits that establish the correspondence are
the points $\{i\}$ and $\{-i\}$, and the real 
projective line $\R\Bbb P^1$.
\end{Example}

\subsection{Statement of results}\label{ssres}
Let $\CK=\C((t))$ be the field of formal Laurent series,
and let $\CO=\C[[t]]$ be the ring of formal power series.
For an algebraic group $G$, we write $G(\CK)$
for its group of $\CK$-valued points, and similarly $G(\CO)$ for its group of $\CO$-valued points.
The quotient set $G(\CK)/G(\CO)$ is known to be the $\BC$-points
of a (not necessarily reduced) ind-finite type complex algebraic ind-scheme. 
In this paper, we shall only be interested in the space of $\BC$-points of this ind-scheme 
equipped with its analytic topology.
We call this space the affine Grassmannian of $G$,
and denote it by $\affgr$.
It is the direct limit of an increasing union of projective complex algebraic varieties.

The aim of this paper 
is to present a version of Theorem~\ref{tmat}
in which the affine Grassmannian 
$\affgr$ plays the role of the 
flag variety $\B$.
The complex algebraic group $K$ is replaced by its group of $\CK$-valued points $K(\CK)$,
and the Lie groups $G_\R$ and $K_c$ are replaced by their
respective polynomial
loop groups $LG_\R$
and $LK_c$. 
These are the subgroups of $G(\BC[t,t^{-1}])$ consisting
of polynomial maps $\C^\times\to G$ which take the unit circle into
the respective Lie group.
Our main result
is the following analogue of the Matsuki correspondence.

\begin{Theorem}\label{taffmat}
For each $LG_\R$-orbit $\orb_\R$ in $\affgr$, 
there exists a unique $K(\ckc)$-orbit $\orb_K$
such that $\orb_\R\cap\orb_K$ is a single $LK_c$-orbit.
This correspondence provides an order-reversing isomorphism
from the poset 
$LG_\R\backslash \affgr$ to the poset $K(\ckc)\backslash \affgr$.
\end{Theorem}
%

\begin{Example}\label{exaffmat}
To simplify the discussion,
instead of considering the lattice of groups of Example~\ref{exmat},
we consider its adjoint form
$$ \begin{matrix} 
 & & \on{SO}_3(\BC) & & \\ 
 & \nearrow & \uparrow & \nwarrow & \\
 \on{SO}_{2,1}(\BR) && \on{SO}_3(\BR) && \on{O}_2(\BC)  \\
 & \nwarrow & \uparrow & \nearrow & \\
 & & \on{O}_2(\BR) & & 
\end{matrix}
$$
and the corresponding lattice of loop groups
$$ \begin{matrix} 
 & & \on{SO}_3(\CK) & & \\ 
 & \nearrow & \uparrow & \nwarrow & \\
L \on{SO}_{2,1}(\BR) && L\on{SO}_3(\BR) && \on{O}_2(\CK)  \\
 & \nwarrow & \uparrow & \nearrow & \\
 & & L\on{O}_2(\BR) & & 
\end{matrix}
$$

The orbit poset $\on{O}_2(\CK) \backslash \affgr$
is isomorphic to the non-negative integers with the standard ordering.
The $\on{O}_2(\CK)$-orbits in $\affgr$ happen to be finite-dimensional
and admit a simple description.
In the group $\on{SO}_3(\BC)$,
consider the maximal torus $\on{SO}_2(\BC)\subset\on{O}_2(\BC)$, and the 
maximal unipotent subgroups
$N^+$ and $N^-$ 
normalized by $\on{SO}_2(\C)$.
Elements of the group $\Hom(\C^\times,\on{SO}_2(\C))$ may be 
canonically identified
with points of $\affgr$.
For each pair of elements $\lambda,\mu\in\Hom(\C^\times,\on{SO}_2(\C))$
thought of as a pair of points of $\affgr$, define 
the subset $\CO^{(\lambda,\mu)}\subset \affgr$ to be the intersection of orbits
$$
\CO^{(\lambda,\mu)}=N^+(\CK)\cdot\lambda\cap N^-(\CK)\cdot\mu.
$$
If $\lambda\neq\mu$,
precisely one of the pair of subsets
$\CO^{(\lambda,\mu)}$ and $\CO^{(\mu,\lambda)}$ is non-empty.
We may fix an identification of the group $\Hom(\C^\times,\on{SO}_2(\C))$ 
with the integers so that $\CO^{(\lambda,\mu)}$ is nonempty if and only if
$\lambda\geq\mu$. Now the $\on{O}_2(\CK)$-orbit $\CO_{K}^n\subset \affgr$
indexed by a non-negative integer $n$ is the disjoint union
$$
\CO_{K}^n=\bigsqcup_{\lambda-\mu=n} \CO^{(\lambda,\mu)}.
$$
When $n$ is zero, the orbit $\CO^n_K$ is isomorphic to the integers, and when $n$ is non-zero,
it is isomorphic to the product of the integers, $\on{SO}_2(\BC)$, and the affine space 
$\BA^{n-1}(\BC)$.

The dual $L\on{SL}_2(\R)$-orbits are infinite-dimensional
and we shall not attempt to describe them here.

The $L\on{O}_2(\BR)$-orbits that establish the correspondence also are finite-dimensional
and admit a simple description.
Fix an 
$\BR$-split maximal torus $T_\BR\subset \on{SO}_{2,1}(\BR)$,
and let $T\subset\on{SO}_3(\BC)$ be its complexification.
As before, elements of the group $\Hom(\C^\times,T)$ may be 
canonically identified
with points of $\affgr$.
For an element $\lambda\in\Hom(\C^\times,T)$ thought of as a point of $\affgr$,
the orbit 
$$
\CO^{\lambda}=L\on{O}_2(\BR)\cdot\lambda
$$
coincides with that for the inverse $-\lambda$. 
After making any identification of the group 
$\Hom(\C^\times,T)$ with the integers,
we obtain an unambiguous $L\on{O}_2(\BR)$-orbit
$$
\CO_{c}^n=\CO^{\lambda}=\CO^{-\lambda}
$$
indexed by the non-negative integer $n$ corresponding to $\lambda$ or $-\lambda$.
These are the $L\on{O}_2(\BR)$-orbits that establish the correspondence.
When $n$ is zero, the orbit $\CO^n_c$ is isomorphic to the integers, and when $n$ is non-zero,
it is isomorphic to the product of the integers and $\on{SO}_2(\BR)$.
\end{Example}

Next we describe the poset $K(\CK)\backslash \affgr$ in general.
Understanding the poset $K\backslash \B$ is a 
difficult problem 
relevant to the representation
theory of $G_\R$~\cite{BB81, MUV92}.
The poset $K(\CK)\backslash \affgr$ turns out to be easy to describe.
To that end, we collect some basic constructions,
and while doing so, establish notation used throughout the paper.

Recall that $\delta$ is a Cartan involution of $G$ which commutes with 
the conjugation $\theta$.
Fix a $\delta$-stable maximal $\R$-split torus $S_\R\subset G_\R$, and
a $\delta$-stable maximal torus $T_{\R}\subset G_\R$ containing $S_\R$.
Let $S$ be the complexification of $S_\R$,
and let $T$ be the complexification of $T_\R$.
Let $\cowt=\morph(\C^\times,T)$ be the lattice
of homomorphisms from $\C^\times$ to $T$,
and let $\rcowt=\morph(\C^\times,S)$ be the sublattice
of homomorphisms with image in $S$.
We call elements of $\cowt$ coweights, and 
elements of $\rcowt$ real coweights. 

Fix a mimimal parabolic subgroup $P_\R\subset G_\R$
containing $T_\R$,
and let $P$ be its complexification.
Let $M_\BR\subset P_\BR$ be the Levi factor of $P_\BR$ containing $T_\BR$,
and let $M\subset P$ be its complexification.
Fix a Borel subgroup $B \subset G$ containing $T$ and contained in $P$.
Let $\dcowt$ be the semigroup of 
coweights 
which are dominant with respect to $B$,
and let $\rdcowt=\rcowt\cap\dcowt$ be the subsemigroup
of real coweights which are dominant.
The lattice $\cowt$, and its subsets, are naturally posets: 
for coweights $\lambda,\mu\in\cowt$, we have $\lambda\leq \mu$ if and only if
$\mu-\lambda$ is a non-negative integral linear combination of positive
coroots of $G$.

Recall that $\eta$ is the involution of $G$ obtained by composing
the conjugation $\theta$ and the Cartan involution $\delta$.
Let $X$ denote 
the identity component
of the subvariety of $G$ on which $\eta$ acts as inverse.
The projection $\pi:G\to X$ defined by $\pi(g)=\eta(g^{-1})g$
allows us to identify $X$ with the symmetric variety $K\backslash G$. 
The torus $S$ lies in $X$ and is maximal among the 
tori of $G$ which lie in $X$.
Consider the diagram 
$$
\pi_1(G)\stackrel{\pi_*}{\to} \pi_1(X) \stackrel{[\cdot]}{\gets} \rdcowt
$$
where the first map is that induced by $\pi:G\to X$,
and the second assigns to a loop
its homotopy class.

\begin{Definition}
Let $\imgp\subset\rdcowt$ be the inverse image via $[\cdot]:\rdcowt\to\pi_1(X)$
of the image of $\pi_*:\pi_1(G)\to\pi_1(X)$.
\end{Definition}

In general, $\imgp$ is a finite index sub-semigroup of $\rdcowt$, and
if $K$ is connected, $\imgp$ is equal to $\rdcowt$.
We order $\imgp$ by restricting the natural order from 
the coweights.

\begin{Theorem}
There is a natural isomorphism of posets 
$K(\ckc)\backslash \affgr\stackrel{\sim}{\to}\imgp$. 
\end{Theorem}

One easily observed feature of the poset $K\backslash \CB$ is the existence
of a unique maximum: if an algebraic group acts on
an irreducible variety with finitely many orbits,
then there exists a dense orbit.
The poset
$K(\CK)\backslash \affgr$ does not contain a maximum unless $K$ contains the derived group of $G$.
In contrast,
it contains a canonical minimum given by
the orbit $K(\CK)\cdot 1G(\CO)\subset\affgr$ corresponding to $0\in \imgp$.
(In general, this is not the unique minimum
since the poset $K(\CK)\backslash \affgr$ is not necessarily connected.)

We conclude the introduction
with some general remarks.
First,
it is worth mentioning that, strictly speaking, the result proved here is not a direct extension
of the Matsuki correspondence formalism to the loop group $G(\BC[t,t^{-1}])$ and the affine Grassmannian $\Gr$.
In general,
the $K(\BC[t,t^{-1}])$-orbits in $\affgr$ are
not discretely parameterized. (See Example~\ref{exaffmat} where the group $K(\BC[t,t^{-1}])$
is equal to the group $K$ itself.)
It is important that we work with the
completion $K(\CK)$.
Second, 
the most enlightening proof of the correspondence of Theorem~\ref{tmat}
involves the Morse theory of moment maps~\cite{MUV92}.
Although we have been motivated by this approach, 
we have not tried to overcome the technical obstacles
in the way of making it work in this setting.
We hope that the proofs given here have the merit of
using other structures which are interesting in their own right. 
Finally, we mention a relation between the orbits of $K(\ck)$ and $LG_\R$ 
in $\affgr$
and the orbits of $K$ and $G_\R$ in a finite-dimensional flag variety of $G$.
For a dominant coweight $\lambda\in\dcowt$, we call the $G$-orbit $C^\lambda=G\cdot\lambda\subset\affgr$
through $\lambda$ the core associated to $\lambda$. The $G$-action on $\gcore$
induces an isomorphism of $\gcore$ with the
flag variety $F(\lambda)=G/P(\lambda)$,
where $P(\lambda)$ is the parabolic subgroup in $G$ which contains
the Borel subgroup $B$, and whose Levi factor is the centralizer of $\lambda$.
By intersecting the orbits of $K(\ck)$ and $LG_\R$ in $\affgr$
with the core $C^\lambda$, we obtain a decomposition of the flag variety $F(\lambda)$.
The orbits of $K$ and $G_\R$ in $F(\lambda)$ refine this decomposition.


\subsection{Overview of paper}\label{sspfs}
The first two sections following the introduction are preparatory.
In Section~\ref{sback}, we recall the necessary
background results concerning the affine Grassmannian
and loop groups.
For example, we recall the fundamental fact that the affine Grassmannian $\affgr$ is homeomorphic
to the based loop group $\Omega G_c$ of based maps $S^1\to G_c$ whose Fourier expansions
are polynomial. 
We also review the natural stratifications of $\affgr\simeq\blg$
and the realizations of the strata as group orbits or via Morse theory.
In Section~\ref{sinv}, we extend 
involutions to loop groups and define the resulting loop subgroups
needed for the Matsuki correspondence.

In Section~\ref{skorbparam},
we parameterize the $K(\CK)$-orbits in $\affgr$.
It is a simple application of the compactification theory of symmetric
varieties. The set of $K(\CK)$-orbits in $\affgr$
is naturally a subset of the set of $G(\CO)$-orbits in $X(\CK)$.
We consider the quotient $G^\ad$
of $G$ by its center $Z(G)$, the corresponding symmetric variety $X^\ad$,
and the wonderful compactification of $X^\ad$.
We are able to parameterize
the $G^\ad(\CO)$-orbits in $X^\ad(\CK)$ in terms of 
the interactions of elements of $X^\ad(\CK)$ with the divisor at infinity
of the wonderful compactification.
Using the projection $X\to X^\ad$, it is then simple to parameterize
the $G(\CO)$-orbits in $X(\CK)$.

In Section~\ref{srorbparam}, we parameterize
the $LG_\BR$-orbits in $\affgr$. It is accomplished by elementary methods.
There is a natural map from the set of $LG_\BR$-orbits in $\affgr$ to the lattice
$\dcowt$. We check directly that the map is injective,
and its image is the subset $\imgp\subset\dcowt$,

The remainder of the paper is devoted to confirming that the order structures
on the orbit posets $K(\CK)\backslash \affgr$ and $LG_\BR\backslash \affgr$ are as asserted.
To accomplish this, we relate the orbits to stratifications
of the based loop space of the symmetric space $K_c\backslash G_c$.
The based loop space of $K_c\backslash G_c$ may be thought of as
honestly based loops in $K_c\backslash G_c$,
or it may be thought of as paths in $G_c$ from the identity to $K_c$.
We consider two models of the based loop space of $K_c\backslash G_c$
which reflect these two possibilities.

First, let $X_c$ denote the identity component
of the subspace of $G_c$ on which the involution $\eta$ acts as inverse. 
The projection $\pi:G_c\to X_c$ defined by $\pi(g)=\eta(g^{-1})g$
allows us to identify $X_c$ with the symmetric space $K_c\backslash G_c$.
In Section~\ref{sbls}, we consider the based loop space $\bls$
of based maps $S^1\to X_c$ whose Fourier expansions
are polynomial. We obtain
stratifications of $\bls$ by restricting those of $\affgr \simeq\blg$.
Under the projection $\pi:\blg\to\bls$, the orbits of 
$K(\CK)$ and $LG_\BR$ in $\affgr \simeq \blg $ are mapped into the strata
of $\bls$.
 In this way, we obtain useful information about the closure relations among the orbits.

In Section~\ref{smatrlg},
we introduce the real affine Grassmannian $\raffgr$. It is 
the real form of $\affgr$ corresponding to the real form $G_\R$ of $G$.
(Although we do not use the result,
it was first recognized by Quillen that $\raffgr$
is homotopy equivalent to the based loop space of the symmetric space $K_c\backslash G_c$.)
As before, we obtain stratifications of $\raffgr$ by restricting those of $\affgr$.
In Sections~\ref{smatbdg} and~\ref{smatspec},
we study a family with general fiber $\affgr$ and special fiber $\raffgr$.
It is a real form of the Beilinson-Drinfeld Grassmannian over a curve.
The family enables us to relate the orbits of $K(\ck)$ and $LG_\R$ in $\affgr$
to the stratifications of $\raffgr$,
and thus to confirm
their closure relations are as asserted.


\subsection{Acknowledgements}
I thank Robert MacPherson for his guidance and support.
I thank Jared Anderson for helpful conversations,
and a referee for many suggestions. I thank Tom Braden for explaining 
the details involved in deducing Theorem~\ref{tcalc}
from known results.
This work was supported in part by the NSF.




\section{Background results}\label{sback}

References for the material in this section include~\cite{Lu83, PS86, BeauLa94, Ginz95}.

\subsection{The affine Grassmannian $\affgr$}\label{ssalg}

Let $G$ be a connected reductive complex algebraic group.
Let $\co=\C[ [ t] ]$ denote the ring of formal power series, 
and let $\ck=\C((t))$ denote its fraction field,
the field of formal Laurent series.
Let $G(\ck)$ denote the group of $\ck$-valued
points of $G$, and let
$G(\co)$ denote the subgroup of $\co$-valued
points.
The quotient set $G(\CK)/G(\CO)$ is naturally the set of $\BC$-points
of a (not necessarily reduced) ind-finite type complex algebraic ind-scheme.
In this paper, we shall only be interested in the space of $\BC$-points of of this ind-scheme equipped
with its analytic topology. We call this space the affine Grassmannian of $G$,
and denote it by $\affgr$. See Remark~\ref{rconcrete} below 
for a concrete realization
of $\affgr$ as the direct limit of an increasing union of projective varieties.


Consider the action of $G(\co)$ on $\affgr$
by left-multiplication. 
Each coweight $\lambda\in\cowt$ determines a point $\lambda\in\affgr$ via
the inclusion $\cowt\subset G(\ck)$.
We call the $G(\co)$-orbits $\gstrat=G(\co)\cdot\lambda\subset\affgr$ 
through the dominant coweights $\lambda\in\dcowt$
the strata of $\affgr$.
For the following, see \cite[Section 11]{Lu83},
or \cite[Theorem (8.6.3) and Proposition (8.6.5)]{PS86}.

\begin{Proposition}\label{pmatstrat}
The affine Grassmannian $\affgr$ is the disjoint union of strata
$$
\affgr=\bigsqcup_{\lambda\in\dcowt} \gstrat.
$$
The stratum $\gstrat$ is a smooth quasi-projective variety. 
Its closure is the disjoint union of strata
$$
\overline\gstrat=\bigsqcup_{\mu\leq\lambda} S^\mu.
$$
\end{Proposition}

\begin{Remark}\label{rconcrete}
The affine Grassmannian
$\affgr$ is the direct limit of its stratum closures
$$
\affgr=\varinjlim_{\lambda} \overline{S^\lambda}. 
$$

\end{Remark}

There is a stratification of 
$\affgr$ dual to the
stratification by $G(\co)$-orbits.
Let $\ipoly$ denote the polynomial ring in the variable $t^{-1}$,
and let $G(\ipoly)$ denote the group of $\ipoly$-valued points of $G$.
Consider 
the action of $G(\ipoly)$ on $\affgr$ by left-multiplication.
We call the $G(\ipoly)$-orbits $\gdstrat=G(\ipoly)\cdot\lambda\subset\affgr$ 
through the dominant coweights $\lambda\in\dcowt$
the dual strata of $\affgr$. 
In general, the dual strata are not finite-dimensional.
For the following, see
\cite[Theorem (8.6.3) and Proposition (8.6.5)]{PS86}.

\begin{Proposition}\label{pmatdstrat}
The affine Grassmannian $\affgr$ is the disjoint union of dual strata
$$
\affgr=\bigsqcup_{\lambda\in\dcowt} \gdstrat.
$$
The closure of the dual stratum $\gdstrat$
is the disjoint union of dual strata
$$
\overline\gdstrat=\bigsqcup_{\lambda\leq \mu} T^\mu.
$$
The intersection of a stratum $S^\mu$ with a dual stratum $\gdstrat$ is nonempty if and only if $\mu\geq \lambda$.

\end{Proposition}

We shall also need the following lemma.
In its proof, and 
throughout this paper, we use the well known fact that the inclusion $G(\lpoly)\hookrightarrow G(\ck)$
induces a bijection 
$$
G(\lpoly)/G(\poly)\risom \affgr
$$ 
(see Theorem~\ref{tfactor} below),
and the group $G(\poly)$ has the same orbits 
in $\affgr$ as the group $G(\co)$.

\begin{Lemma}\label{linter}
The intersection of a stratum $S^\lambda$ with a dual stratum $T^\mu$ is smooth.
\end{Lemma}

\begin{proof}
The quotient set $G(\C((t^{-1})))/G(\C[t])$ is naturally a pro-smooth 
scheme 
which we denote by $\thaffgr$ and call the thick Grassmannian of $G$.
As the inclusion $G(\lpoly)\hookrightarrow G(\ck)$
induces a bijection 
$$
G(\lpoly)/G(\poly)\risom \affgr,
$$ 
the inclusion $G(\lpoly)\hookrightarrow G(\C((t^{-1}))$ induces an
inclusion 
$$
\affgr\hookrightarrow \thaffgr.
$$
As the group $G(\poly)$ has the same orbits 
in $\affgr$ as the group $G(\co)$, to prove the lemma
it suffices to prove that the intersections of the orbits of $G(\C[t])$
and $G(\C[[t^{-1}]])$ in $\thaffgr$
are smooth. 
Any such intersection $Z$ is contained in an open subscheme $Y\subset \thaffgr$
such that there is an integer $k$
for which the kernel $G(t^k\BC[[t]])$ of the projection $G(\BC[[t]])\to G(\BC[[t]]/(t^k))$
acts freely on $Y$, the quotient $G(t^k\BC[[t]])\backslash Y$ is smooth,
and the natural map $Z\to G(t^k\BC[[t]])\backslash Y$ is a closed immersion.
(See \cite[Lemma~2.2.1]{KT95} for a similar result.)
Thus the lemma follows from the transversality
of the subalgebras $\lieg(\C[[t^{-1}]])$ and $\lieg(\C[t])$
in the algebra $\lieg(\C((t^{-1})))$. 
\end{proof}

For each dominant coweight $\lambda\in\dcowt$,
we call the intersection $\gcore=\gstrat\cap \gdstrat$ 
the core associated to $\lambda$.
Let $P(\lambda)$ be the
parabolic subgroup of $G$ which contains the Borel subgroup $B$, and whose Levi factor
is the centralizer of $\lambda$.
The action of $G$ on $\affgr$ preserves the core $\gcore$, 
and $P(\lambda)$ is the stabilizer in $G$ of $\lambda$.

\begin{Proposition}\label{pcore}
For 
$\lambda\in\dcowt$,
the core $\gcore$ 
is isomorphic to the flag variety $G/P(\lambda)$.
\end{Proposition}


\subsection{The based loop group $\blg$}\label{ssgeom}

Let $G_c\subset G$ be a maximal compact subgroup, and
let $S^1\subset \C^\times$ denote the unit circle.
The loop group $LG_c$ is defined to be the subgroup of $G(\lpoly)$
consisting of polynomial maps $\C^\times\to G$ that take $S^1\subset\C^\times$ 
into $G_c\subset G$.
The based loop group $\Omega G_c$ is defined to be the subgroup of $LG_c$
consisting of maps that take $1\in S^1$ to $1\in G_c$.
We endow $\Omega G_c$ with the topology induced by the bijection
of the following theorem.

\begin{Theorem}[\cite{PS86}, Theorem~(8.6.3)]\label{tfactor}
The inclusion $\Omega G_c\subset G(\ck)$
induces a bijection $\Omega G_c\stackrel{\sim}{\to}\affgr$.
\end{Theorem}

The theorem is equivalent to the statement that each element of
$G(\ck)$ factors uniquely into a product of an element of
$\Omega G_c$ and an element of $G(\co)$.
In this paper, we identify the affine Grassmannian $\affgr$ and
the based loop group $\blg$ as topological spaces,
and write $\Omega G_c$ or $\affgr$ depending on the context. 

The choice of a bi-invariant Riemannian metric on $G_c$
provides a natural function on the based loop group $\blg$:
the energy of a loop $\gamma\in\blg$ is defined by the formula
$$
E(\gamma)=\int_{S^1}\|\gamma'\|^2d\theta.
$$ 
The Morse theory of the energy function is
well developed, and
best viewed using the identification
$\blg\simeq\affgr$.
For each $z\in\C^\times$, 
we have the automorphism $t\mapsto z^{-1}t$ of the field $\ck=\C((t))$ of formal Laurent
series.
It induces an automorphism of the group $G(\ck)$ which
descends to an automorphism $\phi_z$ of the affine Grassmannian $\affgr$.

\begin{Theorem}[\cite{PS86}, Theorem~(8.9.9)]\label{tflow}
The cores $\gcore$ are the fixed points of the $\BC^\times$-action on $\affgr$
given by $\phi_z$.
The strata $\gstrat$ are the ascending spaces 
$$
\gstrat=\{x\in\affgr|\lim _{z\to\infty} \phi_z(x)\in \gcore\},
$$
and the dual strata $\gdstrat$ are the descending spaces
$$
\gdstrat=\{x\in\affgr|\lim _{z\to0} \phi_z(x)\in \gcore\}.
$$
\end{Theorem}

The connection with the energy function $E$ is discussed in~\cite[Section~8.9]{PS86}. 
If one is willing to work with a larger class of loops,
one can show that $E$ is indeed a Morse function
on the based loop space of $G_c$.
There is a metric on the based loop space of $G_c$ such that
the gradient flow of $E$
is given by the automorphism $\phi_z$, for $z\in\Rp$.
Our space $\blg$ then consists
of exactly those loops whose upward trajectory exists for all time.
We shall have no need for this general theory, but only for the following.

\begin{Proposition}\label{pintermorse}
For any stratum $S^\lambda$ and dual stratum $T^\mu$,
the restriction of the energy function $E$ to their intersection ${S^\lambda\cap T^\mu}$ has a critical point
if and only if $\lambda=\mu$.
\end{Proposition}


\section{Involutions of loop groups}\label{sinv}

As usual 
let $\theta$ be a conjugation of $G$, and let $\delta$ be a Cartan involution
of $G$
which commutes with $\theta$.
Let $\eta$ denote the composition of $\theta$ and $\delta$.
In this section, we extend these involutions
from $G$ to appropriate loop groups.
All of our loop groups are subgroups of $G(\ck)$,
and so we first describe how to extend involutions to it.

\subsection{Extending algebraic involutions}
The algebraic involution $\eta$ induces an involution 
of the group $G(\cA)$ of $\cA$-valued points of $G$ 
for any $\C$-algebra $\cA$, and  
we denote the induced involution by $\eta$ as well. 
The fixed point subgroup of the induced involution 
is the group $K(\cA)$.
In particular, for the $\BC$-algebra $\CK=\BC((t))$, we have the induced involution $\eta$
of the group $G(\ck)$ with
fixed point subgroup $K(\ck)$.

\subsection{Extending anti-algebraic involutions}
To extend the anti-algebraic involutions $\theta$ and $\delta$ to $G(\cA)$ 
for a $\C$-algebra $\cA$, we need a conjugation $c$ of $\cA$.
In that case, we define involutions $\theta$ and $\delta$ of $G(\cA)$
as follows:
if we think of an element $g\in G(\cA)$ as a map
$$
g:\Spec(\cA) \to G,
$$
then $\theta$ takes it to the composite map
$$
\theta(g):\Spec(\cA) \stackrel{c}{\to} \Spec(\cA) \stackrel{g}{\to}G
   \stackrel{\theta}{\to}G,
$$
and similarly for $\delta$. 
When the $\C$-algebra is $\ckc=\C((t))$,
we take $c$ to be the standard conjugation whose fixed points 
are the field $\ckr=\R((t))$ of real 
formal Laurent series.

For any subgroup of $G(\ck)$ which we shall consider, 
the involutions $\theta$, $\delta$, and $\eta$ restrict to give involutions
of the subgroup. Certain subgroups of $G(\ck)$ have an additional
involution which we describe next.

\subsection{Time-reversal}

The $\C$-algebra $\lpoly$ comes equipped with a natural algebraic
involution $\tau$ which exchanges $t$ and $t^{-1}$.
This induces an involution of the loop group $G(\lpoly)$
which we also denote by $\tau$. 
We write $\thetat$ and $\deltat$ for the involutions of $G(\lpoly)$
obtained by composing $\theta$ and $\delta$ with $\tau$.

\subsection{Subgroup lattices}
To keep track of things, 
we organize the fixed-point subgroups of the involutions
into lattices.
Recall that in the classical setting,
the three commuting involutions $\theta$, $\delta$, and $\eta$ of $G$
provide a lattice of subgroups:
 
$$ \begin{matrix}  
 & & G & & \\ 
 & \nearrow & \uparrow & \nwarrow & \\
 G_\R && G_c && K  \\
 & \nwarrow & \uparrow & \nearrow & \\
 & & K_c & & 
\end{matrix}
$$
Here $G_\R$, $G_c$ and $K$ are the fixed points of $\theta$,
$\delta$ and $\eta$ respectively, 
and $K_c$ is the fixed points of $\theta, \delta$ and $\eta$ simultaneously.

The three involutions $\theta$, $\delta$, and $\eta$ of $G(\ck)$ provide
a lattice of subgroups:

$$ \begin{matrix}  
 & & G(\ck) & & \\ 
 & \nearrow & \uparrow & \nwarrow & \\
 G_\R(\ckr) && G_c(\ckr) && K(\ck)  \\
 & \nwarrow & \uparrow & \nearrow & \\
 & & K_c(\ckr) & & 
\end{matrix}
$$
Here $G_\R(\ckr)$, $G_c(\ckr)$ and $K(\ck)$ are the fixed points of $\theta$,
$\delta$ and $\eta$ respectively, 
and $K_c(\ckr)$ is the fixed points of $\theta, \delta$ and $\eta$ simultaneously.
The groups $G_\R(\ckr)$, $G_c(\ckr)$, and  $K_c(\ckr)$ may be identified
with the groups of $\ckr$-valued points of the respective real algebraic
groups. They are the subgroups of $G(\ck)$ consisting of maps 
$\spec(\ck)\to G$ that take the
punctured formal real line $\Spec(\ckr)\subset\Spec(\ck)$ into the respective real algebraic
groups.

The three involutions $\theta^\tau$, $\delta^\tau$, and $\eta$ of $G(\lpoly)$ provide
a lattice of subgroups:

$$ \begin{matrix}  
 & & G(\lpoly) & & \\ 
 & \nearrow & \uparrow & \nwarrow & \\
 LG_\R && LG_c && K(\lpoly)  \\
 & \nwarrow & \uparrow & \nearrow & \\
 & & LK_c & & 
\end{matrix}
$$
Here $LG_\R$, $LG_c$ and $K(\lpoly)$ are the fixed points of $\thetat$,
$\deltat$ and $\eta$ respectively, 
and $LK_c$ is the fixed points of $\thetat, \deltat$ and $\eta$ simultaneously.
The groups $LG_\R$, $LG_c$ and $LK_c$ are the subgroups of $G(\lpoly)$
consisting of maps $\C^\times\to G$ that take the unit circle $S^1\subset\C^\times$ into
the respective real algebraic groups.

Although we shall have use for many of the groups
in the above two lattices,
it is worth pointing out that the Matsuki correspondence 
for the affine Grassmannian $\affgr$ presented here is a statement about the orbits in $\affgr$
of the groups
in the following mixture of the two lattices:

$$ \begin{matrix} 
 & & G(\ck) & & \\ 
 & \nearrow & \uparrow & \nwarrow & \\
 LG_\R && LG_c && K(\ck)  \\
 & \nwarrow & \uparrow & \nearrow & \\
 & & LK_c & & 
\end{matrix}
$$
For the aims of this paper, this is the lattice which is the
proper analogue of the lattice of the finite-dimensional setting.
In general,
the $K(\BC[t,t^{-1}])$-orbits in $\affgr$ are
not discretely parameterized. (See Example~\ref{exaffmat} where the group $K(\BC[t,t^{-1}])$
is equal to the group $K$ itself.)
It is important that we work with the
completion $K(\CK)$.

\subsection{Anti-involutions}\label{ssinvstrat}
In what follows, the involutions often play a less prominent
role than their compositions with the inverse map. 
We write $\teta$ for the anti-involution of $G(\ck)$
obtained by
composing $\eta$ with the inverse map, and
$\tthetat$ for the anti-involution of $G(\lpoly)$
obtained by
composing $\thetat$ with the inverse map.
It is straightforward to check that when restricted to 
the based loop group $\blg$, the two anti-involutions $\teta$ and $\tthetat$
agree.


\section{$K(\CK)$-orbit parameterization}\label{skorbparam}
Recall that we write $X$ for the identity component
of the subvariety of $G$ on which $\eta$ acts as inverse, and
the projection $\pi:G\to X$ defined by $\pi(g)=\eta(g^{-1})g$
induces a $G$-equivariant isomorphism $K\backslash G\risom X$.
The parameterization of the $K(\CK)$-orbits in $\affgr$ is a simple application
of the compactification theory of symmetric varieties.

\subsection{Wonderful compactification}
We recall here part of the compactification theory of symmetric varieties.
Our reference is~\cite[Section 2]{deC87}.
For the moment, let us assume that the center of $G$ is trivial, or
in other words, that $G$ is a semisimple group of adjoint type.

Consider the Grassmannian $\Gr_d(\fg)$ of $d$-dimensional subspaces in the 
adjoint representation $\fg$
where $d$ is the dimension of $K$. The Lie algebra of $K$
defines a point of $\Gr_d(\fg)$ whose stabilizer in $G$ is precisely $K$.
Therefore we may identify the $G$-orbit through this point 
with the symmetric variety $X$.
The closure $Z\subset \Gr_d(\fg)$ of this orbit is called the wonderful compactification
of $X$.

Recall that the torus $S$ is a subvariety of $X$.
Let $Y\subset Z$ be the partial compactification of $S$
characterized
by the property that for $\lambda\in\rcowt$, we have
$$
\lim_{t\to 0}\lambda(t)\in Y
\mbox{ if and only if }
\lambda\in\rdcowt.
$$
Let $P^u$ be the unipotent radical of the complexification $P\subset G$
of the minimal parabolic subgroup $P_\BR\subset G_\BR$.

\begin{Theorem}[\cite{deC87}, Lemma, p.~372]\label{tlocal}
The action map $P^u\times Y\to Z$ is an open embedding and its image contains
an open nonempty subset of each $G$-orbit in $Z$.
\end{Theorem}

\subsection{Loop parameterization}
We are now ready to parameterize the $K(\CK)$-orbits in $\affgr$.
First we parameterize the $G(\CO)$-orbits in $X(\CK)$.
Note that via the inclusion $S\subset X$,
we may identify $\rdcowt$ with a subset of $X(\CK)$.

\begin{Theorem}\label{tloop}
Each $G(\CO)$-orbit in $X(\CK)$ contains a unique element of $\rdcowt$.
\end{Theorem}

\begin{proof}
First, assume $G$ is semisimple of adjoint type.
Let $Z$ be the wonderful compactification of $X$.
Via the embedding $X\to Z$, we obtain an injection $X(\CK)\to Z(\CK)$.
Since $Z$ is complete, each element $\gamma\in Z(\CK)$ extends
to an element $\bar\gamma\in Z(\CO)$. By Theorem~\ref{tlocal}, 
the $G(\CO)$-orbit through $\bar\gamma$ contains an element $\bar\gamma'\in Z(\CO)$
which lies in the image of the action map $P^u\times Y\to Z$.
Therefore we may consider $\bar\gamma'$ as an element of 
$P^u(\CO)\times Y(\CO)$.
Acting by an element of $P^u(\CO)$,
we see that the $G(\CO)$-orbit through $\gamma$ contains an
element in $Y(\CO)$. 
Acting by an element of $T(\CO)$,
we conclude that the $G(\CO)$-orbit through $\gamma$ contains an
element $\lambda\in\rdcowt$.

Now for $G$ not necessarily semisimple of adjoint type,
we write $G^\ad$ for the quotient of $G$ by its center $Z(G)$,
and $S^\ad$ for the image of $S$ under the projection $G\to G^\ad$.
The involution $\eta$ descends to $G^\ad$, and we write $X^\ad$
for the resulting symmetric space.
Quotienting by the action of $Z(G)$ 
gives a $G$-equivariant map $X\to X^\ad$
which restricts to a $T$-equivariant map $S\to S^\ad$. 
We have proved above that each 
$G^\ad(\CO)$-orbit in $X^\ad(\CK)$ contains an element of $\Lambda_{S^\ad}^+$.
Therefore
each 
$G(\CO)$-orbit in $X(\CK)$ contains an element of $S(\CK)$ which projects
to an element of $\Lambda_{S^\ad}^+$.
Since $\rdcowt$ is the inverse image of $\Lambda_{S^\ad}^+$
under the projection $\rcowt\to\Lambda_{S^\ad}$,
acting by an element of $T(\CO)$,
we conclude that each $G(\CO)$-orbit in $X(\CK)$ contains an
element of $\rdcowt$.

To prove the uniqueness assertion, consider $\gamma\in X(\CK)$ as an element of $G(\CK)$
via the inclusion $X\subset G$.
Since the action of $g\in G$ on $x\in X$ 
is given by $g\cdot x= \eta(g^{-1})x g$,
the orbit $G(\CO)\cdot\gamma\subset X(\CK)$ is contained in 
the double coset $G(\CO)\gamma G(\CO)\subset G(\CK)$.
By Proposition~\ref{pmatstrat}, this double coset contains a unique element of $\dcowt$.
Thus the element of $\rdcowt$ in the $G(\CO)$-orbit
through~$\gamma$ is unique.
\end{proof}

There is a natural inclusion of the $K(\CK)$-orbits in $\affgr$ into the $G(\CO)$-orbits in $X(\CK)$.
To identify its image,
recall the diagram 
$$
\pi_1(G)\stackrel{\pi_*}{\to} \pi_1(X) \stackrel{[\cdot]}{\gets} \rdcowt
$$
where the first map is that induced by $\pi:G\to X$,
and the second assigns to a loop
its homotopy class.
Recall that $\imgp\subset\rdcowt$ is defined to be the inverse image via $[\cdot]:\rdcowt\to\pi_1(X)$
of the image of $\pi_*:\pi_1(G)\to\pi_1(X)$.

\begin{Theorem}
The image of the natural inclusion $K(\CK)\backslash\affgr\to X(\CK)/G(\CO)$
is the subset $\imgp\subset\rdcowt$.
\end{Theorem}

\begin{proof}
Let $\partial:\pi_1(X)\to\pi_0(K)$ denote the boundary map in the long exact sequence in homotopy
of the fibration $\pi:G\to X$. By exactness, to show $[\lambda]\in\pi_1(X)$
is in the image of $\pi_*$ is the same as to show $\partial([\lambda])=[1]$.
Consider the map $\delta:X(K)\to\pi_0(K)$
which assigns to $\gamma\in X(K)$ the class of the pull-back via $\gamma$ of
the $K$-bundle $\pi:G\to X$.
It clearly has the following properties: (i) the composition 
$$
G(\CK)\stackrel{\pi}{\to} X(\CK)
\stackrel{\delta}{\to}\pi_0(K)
$$ 
takes any $g\in G(K)$ to $[1]\in\pi_0(K)$, (ii) for $g\in G(\CO)$,
$\delta(g\cdot \gamma)=\delta(\gamma)$ since $\pi:G\to K$
is $G$-equivariant,
and (iii) the restriction of $\delta$ to $\rdcowt$ agrees with $\partial$.
It follows immediately that
the image of $K(\CK)\backslash\affgr\to X(\CK)/G(\CO)$ is contained in $\imgp$.
The following lemma then implies the theorem.
\end{proof} 

\begin{Lemma}\label{lgeodinimage}
For each $\lambda\in\imgp$, there exists $c\in\blg$ such that $\pi(c)=\lambda$.
\end{Lemma}

\begin{proof}
We would like to take $c(t)$ to be the loop $\lambda(t^{1/2})$.
But if $\lambda(-1)\neq 1$, this will not
be a closed loop. To correct for this, we find
a homomorphism $k:S^1\to K_c$ such that
$k(t^{1/2})\lambda(t^{1/2})$ is in $\blg$.
Let $\partial:\pi_1(X_c)\to\pi_0(K_c)$ denote the boundary map in the long exact sequence in homotopy
of the fibration $\pi:G_c\to X_c$.
By definition, $[\lambda]\in\pi_1(X_c)$ is in the image of $\pi_*$,
and so by exactness, we have $\partial([\lambda])=[1]$.
The path $\lambda(t^{1/2})$ is a lift with which
we may calculate $\partial$, so $\lambda(-1)$
belongs to the identity component of $K_c$.
Therefore since $\lambda(-1)^2=1$, we may define a homomorphism
$k:S^1\to K_c$ such that
$k(-1)=\lambda(-1)$.

We check that $c(t)=k(t^{1/2})\lambda(t^{1/2}) $ is in $\blg$.
First realize $G_c$ as a matrix group so that $\lambda(t)$ is a diagonal
matrix whose diagonal entries are powers of $t$. 
Hence $\lambda(-1)$ is a diagonal matrix 
with each diagonal entry equal to $+1$ or $-1$.  
Now observe that since $k(t)$ is a homomorphism
and $k(-1)=\lambda(-1)$,
its image must lie in the centralizer $L$ of $\lambda(-1)$ for all $t$.
Each entry of $L$ is a block diagonal matrix with
two blocks, $B_+$ for the $+1$ block of $\lambda(-1)$ and
$B_-$ for the $-1$ block. 
The assertion that $c(t)$ is in $\blg$
is equivalent to the assertion (${\dag}$): 
the entries of $k(t^{1/2})$ in $B_+$ are in $\lpoly$ 
and the entries in $B_-$ are in $t^{1/2}\lpoly$.
Choose $\ell\in L$ such that $\ell^{-1}k(t^{1/2})\ell$ is diagonal for all $t$.
Since $k(-1)=\lambda(-1)$, the path $\ell^{-1}k(t^{1/2})\ell$ satisfies ($\dagger$).
Since the entries of $\ell$ are independent of $t$,
we conclude that $k(t^{1/2})$ satisfies ($\dagger$) and the lemma is proved.
\end{proof}


\section{$LG_\BR$-orbit parameterization}\label{srorbparam}

The aim of this section is to show that the subset $\imgp\subset \rdcowt$ also parameterizes the
$LG_\BR$-orbits in $\affgr$.
We begin by constructing a natural map from the set of $LG_\BR$-orbits in $\affgr$ to the set $\dcowt$.
It will then take some work to show that the map is an inclusion with image $\imgp$.
Recall that we write $\tthetat$ for the anti-involution of $G(\BC[t,t^{-1}])$
obtained by composing $\thetat$ with the inverse map.

Given an element $g\in G(\BC[t,t^{-1}])$, consider the element $\ttheta^\tau(g^{})g\in G(\BC[t,t^{-1}])$.
By Proposition~\ref{pmatdstrat}, the double coset $G(\BC[t^{-1}])\ttheta^\tau(g^{})g G(\BC[t])$
contains a unique element $\lambda\in\dcowt$.
Since $\theta^\tau$ fixes elements of $LG_\BR$ and takes $G(\BC[t])$ to $G(\BC[t^{-1}])$,
the element $\lambda$ only depends on the orbit of $g$ under the action of $LG_\BR$ by
left multiplication and $G(\BC[t])$ by right multiplication.
Therefore $\lambda$ is an invariant of the $LG_\BR$-orbit through the point in $\affgr$
represented by $g$. 
In other words, to each $LG_\BR$-orbit in $\affgr$, we have associated an element of $\dcowt$.

\begin{Proposition}\label{prparam}
The image of the map
$
LG_\BR\backslash \affgr\to\dcowt
$
is contained in the subset $\rdcowt\subset\dcowt.$
\end{Proposition}

\begin{proof}
We first collect some results on how the anti-involution
$\tilde\theta^{\tau}$ acts on the coweight lattice $\cowt$.
Note that the involution $\tau$ and the inverse map coincide when restricted
to $\cowt$, and so $\tilde\theta^{\tau}$ coincides with the involution induced
by the conjugation $\theta$.
Therefore the fixed points of $\tilde\theta^{\tau}$ are precisely the real coweights $\rcowt\subset\cowt$.

The action of $\tilde\theta^{\tau}$ on $\cowt$ induces an action on the dominant coweights
$\dcowt$. By definition, it takes a dominant coweight $\lambda$ to the dominant coweight 
in the Weyl group orbit through $\tilde\theta^\tau(\lambda)$.
Although the fixed points of $\tilde\theta^{\tau}$
acting on $\dcowt$
contain the real dominant coweights $\rdcowt$, they also in general contain other dominant coweights.
Indeed, the action of $\tilde\theta^{\tau}$ on $\dcowt$ coincides with that
of the involution $\tilde\kappa^\tau$ for the quasi-split conjugation $\kappa$
in the same inner class as $\theta$. More precisely,
the two involutions acting on $\cowt$ are related by the formula
$$
\tilde\theta^\tau=\tilde\kappa^\tau\circ w_M
$$
where $w_M$ is the longest element of the Weyl group of the complexification $M\subset P$
of the Levi factor $M_\BR\subset P_\BR$ of the minimal parabolic subgroup
$P_\BR\subset G_\BR$.
Thus an element $\lambda\in\dcowt$ is in fact in $\rdcowt$ if and only if 
(i) as an element of $\dcowt$, it
is fixed by $\tilde\theta^\tau$, and (ii) as an element of $\cowt$, it is fixed by
$w_M$.

Now fix $x\in\affgr$ with representative $gG(\BC[t])$,
 and let $\lambda\in\dcowt$ be contained in 
$G(\BC[t^{-1}])\ttheta^\tau(g^{})g G(\BC[t])$.
Clearly $G(\BC[t^{-1}])\ttheta^\tau(g^{})g G(\BC[t])$ is preserved
by $\tilde\theta^{\tau}$, and so $\lambda$ and $\tilde\theta^{\tau}(\lambda)$
are conjugate under the Weyl group of $G$. 
Thus by definition, $\lambda$ 
is fixed by the action of $\tilde\theta^\tau$ on $\dcowt$.
To prove the proposition, it remains to show that $\lambda$ 
is fixed by the action of $w_M$ on $\cowt$.

Recall that to $\lambda\in\dcowt$, we associate the parabolic subgroup $P(\lambda)\subset G$
which contains the Borel subgroup $B$, and
whose Levi factor is the centralizer of $\lambda$.
Observe that $\lambda$ is fixed by $w_M$ if and only if
$P(\lambda)$ 
contains the parabolic subgroup $P$.
To show that $P(\lambda)$ contains $P$, 
it suffices to find a parabolic subgroup $Q\subset G$ conjugate to $P(\lambda)$
and preserved by the conjugation $\theta$.

Consider the``opposite" affine Grassmannian $\affgr_-=G(\BC[t,t^{-1}])/G(\BC[t^{-1}])$.
By Proposition~\ref{pmatdstrat},
the orbits of the diagonal action of $G(\BC[t,t^{-1}])$ on the product $\affgr\times\affgr_-$
are indexed by $\dcowt$.
To be precise, given a pair $(x,y)\in\affgr\times\affgr_-$ with representative $(gG(\BC[t]),hG(\BC[t^{-1}]))$, 
the diagonal $G(\BC[t,t^{-1}])$-orbit through
it is indexed by the unique $\lambda\in\dcowt$ contained in the double coset
$G(\BC[t^{-1}])h^{-1}gG(\BC[t])$.
In this case, we say that $\lambda$ is the relative position of the pair $(x,y)$.
In particular, the relative position of the pair $(x,\theta^\tau(x))$
is the $\lambda\in\dcowt$ indexing the $LG_\BR$-orbit through $x\in\affgr$.

For a pair $(x,y)\in\affgr\times\affgr_-$, 
let $St(x)\subset G(\BC[t,t^{-1}])$ be the stabilizer of $x\in\affgr$,
and let $St^-(y)\subset G(\BC[t,t^{-1}])$ be the stabilizer of $y\in\affgr_-$.

\begin{Lemma}
Fix a pair $(x,y)\in\affgr\times\affgr_-$, and 
let $\lambda$ be the relative position of $(x,y)$.
For any $a\in \BC^\times$, 
the subgroup 
$$
Q=\{g\in G | g=s(a), \mbox{ for } s\in St(x)\cap St^-(y)\}\subset G
$$
is conjugate in $G$ to the parabolic subgroup $P(\lambda)$.
\end{Lemma}

\begin{proof}
Under the diagonal action of $g\in G(\BC[t,t^{-1}])$
on the product $\affgr\times\affgr_-$, the subgroup $Q$ is conjugated by $g(a)$. 
Therefore by Proposition~\ref{pmatdstrat}, we may assume $x$
has representative $\lambda G(\BC[t])$ and $y$ has representative $1G(\BC[t^{-1}])$.
For such a pair, it is easy to verify that $Q$ is precisely $P(\lambda)$.
\end{proof}

Now we apply the above to the pair $(x,\theta^\tau(x))$. The intersection
$St(x)\cap St^-({\theta^\tau(x)})$ is preserved by $\theta^\tau$, and for any $a\in S^1$,
evaluating at $a$ is an equivariant map with respect to the action of $\theta^\tau$
on $G(\BC[t,t^{-1}])$ and $\theta$ on $G$. Therefore the subgroup
$$
Q=\{g\in G | g=s(a), \mbox{ for } s\in St(x)\cap St^-(\theta^\tau(x))\}\subset G
$$
is a parabolic subgroup conjugate to $P(\lambda)$ and preserved by $\theta$.
This completes the proof of the proposition.
\end{proof}

We are now ready to parameterize the $LG_\BR$-orbits in $\affgr$.

\begin{Theorem}
The map $LG_\BR\backslash \affgr\to\dcowt$ is an inclusion.  Its image is the subset $\imgp\subset\dcowt$.
\end{Theorem}

\begin{proof}
Fix $x\in\Omega G_c$,
and let $y\in \Omega G_c$
be the element $\ttheta^\tau(x)x$.
By Proposition~\ref{prparam}, 
the unique $\lambda\in\dcowt$ contained in the double coset
$G(\BC[t^{-1}]) yG(\BC[t])\subset G(\BC[t,t^{-1}])$ is in fact in $\rdcowt$.

We first show that the following assertion ({\ddag}) 
implies 
the theorem: 
$$
\mbox{
There exists $g_+\in G(\poly)$ such that $\tthetat(g_+)yg_+=\lambda$ in $G(\lpoly)$.
}
$$
Let $X_\BR\subset G$ denote the identity component of the fixed-point
subspace of $\ttheta$. Thinking of $y$ and $\lambda$ as loops in $X_\BR$, we see that
the assertion ({\ddag}) implies that $[y]=[\lambda]$ in $\pi_1(X_\BR)$ since 
$g_+$ extends across $0$ and $\tthetat(g_+)$ extends across $\infty$.
Thinking of $y$ and $\lambda$ as loops in $X$,
we conclude that $[y]=[\lambda]$ in $\pi_1(X)$
since $X_c$ is a deformation retract of both $X$ and $X_\BR$.
This proves that $\lambda\in\imgp$,
since $\ttheta^\tau$ and $\teta$ agree on $\Omega G_c$,
and so by construction $[y]$ is in the image of $\pi_*$.
Now by Lemma~\ref{lgeodinimage} and the fact that $\ttheta^\tau$ and $\teta$ agree on $\Omega G_c$,
we may choose $c\in\Omega G_c$ such that $\ttheta^\tau(c)c=\lambda$.
Since $y=\ttheta^\tau(x)x$ and $\lambda=\ttheta^\tau(c)c$ in $G(\BC[t,t^{-1}])$,
the assertion ({\ddag}) implies that $xg_+c^{-1}$ is fixed by the involution $\theta^\tau$.
In other words, the element $g_\BR=xg_+c^{-1}$ is in $LG_\BR$, we have $g_\BR c=xg_+$
in $G(\BC[t,t^{-1}])$, and so $g_\BR\cdot c=x$ in $\affgr$.
Since the subgroup $\Omega K_c\subset LG_\BR$ acts transitively on the elements
$c\in\Omega G_c$ with $\tthetat(c)c=\lambda$, we conclude that the $LG_\BR$-orbit
in $\affgr$ through $x$ is the only $LG_\BR$-orbit indexed by $\lambda$.
This proves the theorem once we establish the assertion (\ddag).

We shall establish the assertion ({\ddag}) 
by acting on $y$ by a series of elements in $G(\BC[t])$
until we arrive at $\lambda$. 

First, we write $y$ as a product
$g_1\lambda g_2$ with $g_1\in G(\ipoly)$ and $g_2\in G(\poly)$.
By acting by $g_1^{-1}$ on the left and $\tthetat(g_1^{-1})$ on the right, 
we may assume $y=\lambda g_3 $ with $g_3\in G(\poly)$.

Let $L\subset G$ be the centralizer of $\lambda$.

\begin{Lemma}
There is a parabolic subgroup $Q\subset G$ with Levi factor $L$
such that $g_3=\ell u$ with 
$\ell \in L$ and $u \in Q^u(\poly)$, where $Q^u$ is
the unipotent radical of $Q$.
Furthermore, $\ell$ is fixed by $\ttheta$.
\end{Lemma}

\begin{proof}
Applying $\ttheta^\tau$ to $y$
and using that $y$ and $\lambda$ are fixed by $\ttheta^\tau$, we obtain the identity 
$$
\ttheta^\tau(g_3)=\lambda g_3\lambda^{-1} \mbox{ in } G(\BC[t,t^{-1}]).
$$
From the identity and the fact that $g_3\in G(\poly)$ and $\ttheta^\tau(g_3)\in G(\ipoly)$,  
we have the first assertion.
We see the second by inspecting 
the constant term of the identity. 
\end{proof}

Let $g_4=\ell^{-1} g_3$, so that $g_4\in Q^u(\ipoly)$.

\begin{Lemma}
We may define an element $g_4^{1/2}\in G(\ipoly)$ by the formula
$$
g_4^{1/2}=\exp(\log(g_4)/2).
$$
It satisfies the equation
$$
\tthetat(g_4^{-1/2})y g_4^{-1/2}=\lambda\ell. 
$$
\end{Lemma}

\begin{proof}
We may define $g_4^{1/2}$ since $g_4\in Q^u(\ipoly)$.
The equation follows from a simple calculation and 
the observation
that commuting relations among group elements are satisfied by power series
in those elements.
\end{proof}

After acting by $\tthetat(g_4^{-1/2})$ on the left and $g_4^{-1/2}$ on the right,
we may assume $y=\lambda\ell$, with $\ell\in L$ 
satisfying $\ttheta(\ell)=\ell$.
We are guaranteed that $\ell\in X_\BR$ since our original element $y$ was a loop in $X_\BR$.
Since we may identify $X_\BR$ with the quotient space $G_\BR\backslash G$,
we see 
that there exists $h\in G$ such that $\ttheta(h^{})\ell h=1$. 
After acting by $\ttheta(h^{})$ on the left and $h$ on the right,
we may assume $y=\ttheta(h^{})\lambda \theta(h)$
which is a homomorphism $\BC^\times\to G$.

On homomorphisms $\BC^\times\to G$, the time involution $\tau$ and the inverse map agree.
Therefore since $y$ is fixed by $\tthetat$, its values on $\BR^\times$ are fixed by $\theta$.
In other words, $y$ is the extension of a homomorphism $\BR^\times\to G_\BR$.
Therefore it is conjugate to $\lambda$ by an element $g_\BR\in G_\BR$,
and we have $\ttheta(g_\BR)y g_\BR=g_\BR^{-1} y g_\BR=\lambda$,
and the assertion ({\ddag}) is proved.
\end{proof}


\section{The based loop space $\Omega X_c$}\label{sbls}
Recall that we write $X_c\subset G_c$ for the identity component
of the fixed-point subspace of the anti-involution $\teta$, and
the projection $\pi:G_c\to X_c$ defined by $\pi(g)=\teta(g)g$
induces a $G_c$-equivariant isomorphism $K_c\backslash G_c\risom X_c$.

\subsection{Definition of $\bls$}
We define the loop space $LX_c$ to be the subspace of $LG_c$
consisting of maps that take the unit circle $S^1\subset\C^\times$
into $X_c\subset G_c$.
We define the based loop space $\Omega X_c$ to be the subspace of $LX_c$
consisting of maps that take $1\in S^1$ to $1\in X_c\subset G_c$.

\begin{Lemma}
The based loop space $\bls\subset\blg$ is the fixed-point subspace 
of the anti-involution $\teta$.
\end{Lemma}

\begin{proof}
Clearly every element of $\bls$ is fixed by $\teta$. 
Conversely, suppose an element of $\blg$ is fixed by $\teta$. 
Since the element is based,
it maps $S^1$ into the identity component of
the fixed-point subspace of $\teta$.
Thus by definition the element lies in $\bls$.
\end{proof}


The point of introducing $\Omega X_c$ is that 
the $G_c$-equivariant projection
$
\pi: G_c\to X_c
$
induces an $LG_c$-equivariant projection
$
\pi:\Omega G_c\to\Omega X_c
$
which in turn induces an $LG_c$-equivariant inclusion
$
\Omega K_c\backslash \Omega G_c\subset \Omega X_c.
$
By definition, we have an $LG_c$-equivariant inclusion $\Omega X_c\subset\Omega G_c$.
Therefore we may relate structures on $\Omega G_c$
to those on $\Omega K_c\backslash \Omega G_c$.

\subsection{Application to orbit posets}
We have shown that the orbits of $K(\CK)$
and $LG_\BR$ in $\affgr$ are parameterized by the set $\imgp\subset\rdcowt$.
For $\lambda\in\imgp$,
we write $\korb\subset\affgr$ for the $K(\CK)$-orbit indexed by $\lambda$,
and $\rorb\subset\affgr$ for the $LG_\BR$-orbit indexed by $\lambda$.
We also write $\corb\subset\affgr$ for their intersection $\korb\cap\rorb$.

Similarly, for $\lambda\in\imgp$,
we define the subspace $\xstrat\subset\bls$
to be the intersection of $\bls$ with the stratum $\gstrat\subset\blg$,
and the dual subspace $\xdstrat\subset\bls$ 
to be the intersection of $\bls$ with the dual stratum $\gdstrat\subset\blg$.
We also define the core $\xcore\subset\bls$ 
to be the intersection of $\bls$ with the core $\gcore\subset\blg$.

\begin{Proposition}\label{porbinto}
Under the projection 
$\pi:\blg\to\bls,$ 
the $K(\CK)$-orbit $\korb$ maps into the subspace $\xstrat$,
the $LG_\BR$-orbit $\rorb$ maps into the dual subspace $\xdstrat$,
and their intersection $\corb$ maps into the core $\xcore$.
\end{Proposition}

\begin{proof}
By construction, for a point $g\in\korb$ thought of as an element of $\blg$, 
we have that $\pi(g)=\teta(g)g$ is in the double coset
$G(\poly)\lambda G(\poly)\subset G(\BC[t,t^{-1}])$. By definition,
$\xstrat$ is the intersection of $\bls$ with this double coset.

The proof for the $LG_\R$-orbits is similar, keeping in mind 
that the restrictions of 
$\teta$ and $\tthetat$ to $\blg$ are equal.
The assertion for the orbit intersections is then immediate from the definitions.
 \end{proof}
%


We have the following application to the orbit poset structures.

\begin{Proposition}\label{pposetbounds}
For $\lambda,\mu\in\imgp$, we have
$$
\CO_K^\mu \cap\overline \korb=\emptyset \quad\mbox{ if $\mu\not\leq\lambda$}
$$
$$
\overline{\rorb}  \cap \CO_\BR^\mu =\emptyset \quad\mbox{ if $\lambda\not\leq\mu$.}
$$
\end{Proposition}

\begin{proof}
From the definitions and Propositions~\ref{pmatstrat}~and~\ref{pmatdstrat}, we have
$$
P^\mu \cap\overline \xstrat=\emptyset \quad\mbox{ if $\mu\not\leq\lambda$}
$$
$$
\overline{\xdstrat}  \cap Q^\mu =\emptyset \quad\mbox{ if $\lambda\not\leq\mu$.}
$$
The result then follows immediately from Proposition~\ref{porbinto}.
\end{proof}

The following is also useful.

\begin{Proposition}\label{pstratman}
For $\lambda,\mu\in\imgp$,
the subspace $\xstrat\subset\bls$ and its intersection $\xstrat\cap Q^\mu$ with the dual subspace 
$Q^\mu \subset\bls$ are smooth.
The restriction of the energy function $E$ to the intersection ${P^\lambda\cap Q^\mu}$
has no critical point if $\lambda\not =\mu$. 
\end{Proposition}

\begin{proof}
We use the following easily verified lemma.

\begin{Lemma}\label{leqmorse}
Let $N$ be a smooth manifold on which a compact group $C$
acts smoothly, and let $M\subset N$ denote the fixed-point subspace of $C$. 
Then $M$ is a smooth submanifold.
Let $f:N\to\R$ be a smooth $C$-invariant function,
and let $f|_M$ denote the restriction of $f$ to $M$.
If $m\in M$ is a critical point of $f|_M$, then $m$ is also a critical
point of $f$. 
\end{Lemma}

That $\xstrat$ is smooth follows immediately from Proposition~\ref{pmatstrat} and 
Lemma~\ref{leqmorse} with $N=\gstrat$,
$M=\xstrat$, and $C=\langle\teta\rangle$.
The remaining assertions follow immediately from Lemma~\ref{linter}, Proposition~\ref{pintermorse},
and  
Lemma~\ref{leqmorse} with $N=S^\lambda\cap T^\mu$,
$M=P^\lambda\cap Q^\mu$, $C=\langle\teta\rangle$, and $f=E|_{S^\lambda\cap T^\mu}$.
\end{proof}

As a corollary of the above propositions, we obtain
the geometric relation between the orbits of $K(\CK)$
and $LG_\BR$ in $\affgr$ which characterizes the correspondence.

\begin{Corollary}\label{cmatcorr}
For $\lambda,\mu\in\imgp$, the intersection $\korb\cap \orb_\R^\mu$ of 
the $K(\ckc)$-orbit $\korb\subset\affgr$
and the $LG_\R$-orbit $\orb_\R^\mu\subset\affgr$, 
is a single $LK_c$-orbit
if and only if $\lambda=\mu$.
\end{Corollary}

\begin{proof}
By Proposition~\ref{porbinto} and the fact that $K_c$ acts transitively
on the core $\xcore$, the inverse image $\CO_c^\lambda=\pi^{-1}(\xcore)$ is
a single $LK_c$-orbit and the unique one contained in $\rorb\cap\korb$.
On the other hand, 
if $\lambda\neq\mu$ and the intersection $\xstrat\cap Q^\mu$ is nonempty, 
then it is not compact by Proposition~\ref{pstratman}.
Since $\pi$ maps each $LK_c$-orbit in $\blg$ to a single $K_c$-orbit 
in $\bls$, there is an infinite number of $LK_c$-orbits in $\orb^\lambda_K\cap\orb^\mu_R$.
\end{proof}

In the remainder of this section, our aim is to prove a weak converse to Proposition~\ref{pposetbounds}.
To do this, we first need the following.

\begin{Proposition}\label{pcomponto}
For $\lambda\in\imgp$, let $P^\lambda_0$ be a connected component of the subspace $\xstrat\subset\bls$ such that a point of $P^\lambda_0$ is in the image of the projection $\pi:\blg\to\bls$.
Then there is a compact subspace $W\subset\blg$ such that 
$P^\lambda_0\subset\pi(W)$.
\end{Proposition}

\begin{proof}
We use the following easily-verified lemma.

\begin{Lemma}\label{leqsurj}
Let $N$ 
be a smooth manifold on which a compact group $C$ acts smoothly
with fixed point set $M\subset N$.
Let $Z$ be a topological space on which $C$ acts continuously
with fixed point set $Y\subset Z$.
Assume $Z$ is the disjoint union of smooth manifolds $Z=\cup_{i\in I} Z_i$,
such that the action of $C$ preserves $Z_i$ and is smooth, for each $i\in I$.
Let $f:Z\to N$ be a proper $C$-equivariant map such that
for each $i\in I$, the restriction $f|_{Z_i}$
is a smooth submersion.
Then the restriction $f|_Y$ is a surjection onto the components of $M$
which contain a point in the image of $f|_Y$.

\end{Lemma}

We shall apply Lemma~\ref{leqsurj} with $N=\gstrat$, $M=P^\lambda$, and $C=\langle\teta\rangle$. 
We must construct the other ingredients required by the
lemma.

Consider the product $\blg\times\blg$ 
of two copies of the based loop group.
Let $m$ denote the multiplication map
$$ \begin{matrix}  
  m:\blg\times\blg\to\blg, & m((x,y))=xy.
\end{matrix}
$$
Let $p$ denote the convolution map
$$ \begin{matrix}  
  p:G(\ckc)\times_{G(\coc)}{G(\ckc)/G(\coc)}\to \affgr, &
  p(g,hG(\coc))=ghG(\coc)
\end{matrix}
$$
where $G(\ck)\times_{G(\co)}{G(\ck)/G(\co)}$ is the twisted
product defined by the equivalence relation 
$(xg_+,yG(\coc))=(x,g_+yG(\coc))$, for $g_+\in G(\coc)$.
By Theorem~\ref{tfactor}, 
the map $\blg\times\blg\to G(\ck)\times_{G(\co)}{G(\ck)/G(\co)}$
induced by the inclusion $\blg\times\blg\subset G(\ck)\times G(\ck)$
is a homeomorphism. 
We obtain a commutative diagram
$$ \begin{matrix}  
  \blg\times\blg & \stackrel{\sim}{\to} & G(\ckc)\times_{G(\coc)}{G(\ckc)/G(\coc)} \\ 
  \\
  m\downarrow & & \downarrow p \\
  \\
  \blg & \stackrel{\sim}{\to} & \affgr.
\end{matrix}
$$

The reason for introducing the twisted product is 
that the group $G(\co)$ acts on it by the formula
$g_+\cdot(x,y)=(g_+x,y)$, and the convolution map $p$ is equivariant
for the action.
 
\begin{Lemma}\label{leqsub}
For 
$(x,y)\in G(\ck)\times_{G(\co)}{G(\ck)/G(\co)}$, 
let $S$ denote the $G(\co)$-orbit through 
$(x,y)$,
and let $\gstrat$ denote the stratum of $\affgr$ containing 
$p((x,y))$.
Then the restriction of the convolution $p|_S:S\to \gstrat$ 
is a surjective submersion.
\end{Lemma}

\begin{proof}
This follows immediately from the $G(\co)$-equivariance 
of the convolution.
\end{proof}

To apply Lemma~\ref{leqsurj}, we need a candidate for the space $Z$.
By Lemma~\ref{lgeodinimage}, 
we may find a point $(\teta(x),x)\in G(\ck)\times_{G(\co)}{G(\ck)/G(\co)}$
such that $p((\teta(x),x))=\lambda$.
Consider the $G(\co)$-orbit $S\subset G(\ck)\times_{G(\co)}{G(\ck)/G(\co)}$ 
through the point $(\teta(x),x)$.
Let $Z$ be 
the intersection of the closure $\overline S$ 
and the inverse image $p^{-1}(\gstrat)$.
Clearly it is a disjoint union of $G(\co)$-orbits,
and the map $p:Z\to\gstrat$ is proper and restricts to 
a submersion on each orbit by Lemma~\ref{leqsub}.

Finally, to apply Lemma~\ref{leqsurj},
we are left to find an action of $C=\z/2\z$
on $Z$ such that the convolution $p$
is equivariant for it and 
the action of $\teta$ on $\gstrat$.
Let $\tepsilon$ be the anti-involution of $\blg\times\blg$ defined by
$\tepsilon(x,y)=(\teta(y),\teta(x))$. 
Clearly the multiplication $m$ is equivariant for the action of
$\tepsilon$ on $\blg\times\blg$ and the action of $\teta$ on $\blg$.
Using the identifications $\blg\times\blg\simeq G(\ck)\times_{G(\co)}{G(\ck)/G(\co)}$
and $\blg\simeq\affgr$,
we see that the convolution $p$ is equivariant for
the induced actions.

\begin{Lemma}\label{lepsorb}
For $(x,y)\in G(\ck)\times_{G(\co)}{G(\ck)/G(\co)}$, the map $\tepsilon$ 
takes the $G(\co)$-orbit through
$(x,y)$ to the $G(\co)$-orbit through 
$(\teta(y),\teta(x))$.
\end{Lemma}

\begin{proof}
Let $g_+\in G(\co)$, and let 
$(a,b)=g_+\cdot (x,y)\in G(\ck)\times_{G(\co)}{G(\ck)/G(\co)}$.
We need to show that there exists $h_+\in G(\co)$ such that 
$h_+\cdot (\teta(y),\teta(x))= (\teta(b),\teta(a))$.
By Thereom~\ref{tfactor}, there exists $g'_+\in G(\co)$
such that $g_+xy=abg'_+$ in $G(\ckc)$.
A simple calculation shows that taking $h_+=\teta(g'_+)^{-1}$ works.
\end{proof}

Using Lemma~\ref{lepsorb} for each $G(\co)$-orbit in $Z$, 
we see that $\tepsilon$ acts on $Z$ 
satisfying the requirements of Lemma~\ref{leqsurj}.
We conclude by Lemma~\ref{leqsurj} that the fixed points of $\tepsilon$
in $Z$ surject onto the components $P^\lambda_0$ containing a point in the image
of the convolution $p$.
Now we may identify the fixed points of $\tepsilon$ 
in $G(\ck)\times_{G(\co)}{G(\ck)/G(\co)}$ with $\blg$
so that $p$ restricted to $\blg$
coincides with the projection $\pi$.
By construction, the closure $\overline Z\subset\blg$ is compact,
and we may take $W$ to be $\overline Z$. This completes the proof of the proposition.
\end{proof}

We are now ready to prove a weak converse to Proposition~\ref{pposetbounds}.
We call elements $\lambda, \mu\in\imgp$
primitively related, and write $\lambda\prim \mu$,
if $\lambda\leq \mu$ and there is no element
$\nu\in\imgp$ with $\lambda<\nu<\mu$.

\begin{Proposition} \label{ppriminter}
Let $\lambda\prim \mu\in\imgp$ be primitively related. 
Assume the intersection $\CO_K^\mu\cap \rorb$ is nonempty. Then
{$$
\overline{\CO_K^\mu} \cap\korb\not=
\emptyset 
$$
$$
\overline{\rorb}  \cap \CO_\BR^\mu \not=
\emptyset.
$$}
\end{Proposition}

\begin{proof}
By assumption and Proposition~\ref{porbinto}, 
there is a nonempty component $N$ of the intersection $P^{\mu}\cap \xdstrat$.

We show that $N$
contains an open line
segment whose lower limit lies in the core $\xcore$ and
whose upper limit lies in the core $B^\mu$. 
By Proposition~\ref{pstratman},  $N$ is a manifold,
and the restriction of the energy function $E|_N$
has no critical points.
Thus to find such a line segment, it suffices to show that for $a,b\in\R$
with $E(\lambda)<a\leq b<E(\mu)$, the inverse image 
$(E|_{N})^{-1}([a,b])$ is compact.
If this is so, then for any metric on $N$,
any integral line of the gradient vector field of 
$E|_{N}$ will serve as the sought-after line segment.

We show that the inverse image
$(E|_{N})^{-1}([a,b])$ is compact by contradiction.
Suppose there is a sequence of points 
in it such that no subsequence
has a limit. Then by Propositions~\ref{pmatstrat} and~\ref{pmatdstrat}, 
there is a dominant
coweight $\nu\in\dcowt$ with $\lambda<\nu<\mu$ such that the sequence has a subsequence
with a limit in the stratum $S^\nu$ of $\blg$. 
By Proposition~\ref{pcomponto}, there is a compact
subspace $W\subset \blg$ such that $N\subset \pi(W)$.
Therefore $\pi(W)$ 
contains the limit point,
and so by Proposition~\ref{porbinto}, we have $\nu\in\imgp$.
But this contradicts the assumption
that $\lambda\prim \mu\in\imgp$ are primitively related.

Finally, the line segment we have constructed in $N$ is contained in $\pi(W)$.
Consider any set-theoretic lift of the line segment
to $W$. By the compactness of $W$, in each direction of the lift there is a subsequence whose
limit exists. By construction, one limit projects to $\xstrat$ and the other to $Q^\mu$.
Therefore by Proposition~\ref{porbinto}, one limit lies in $\korb$ and the other in $\CO_\BR^\mu$.
\end{proof}


\section{The real affine Grassmannian $\raffgr$}\label{smatrlg}

The aim of the remaining sections
is to establish that the orbit posets of $LG_\R$ and $K(\ck)$
acting on $\affgr$ are as asserted in the introduction. 
To accomplish this, it suffices by Proposition~\ref{ppriminter} to show that for 
$\lambda\prim \mu\in\imgp$ primitively related,
the intersection $\CO_K^\mu\cap \rorb$ is nonempty.
Proving this intersection is nonempty is the aim of the remaining sections.
The strategy is to relate
the orbits of $LG_\R$ and $K(\ck)$
in $\affgr$ to stratifications of the real affine Grassmannian $\raffgr$
which we introduce in this section.

\subsection{Definition of $\raffgr$}
As usual, let $c$ denote the standard conjugation of the field 
$\ck=\C((t))$ of
formal Laurent series,
and let $\ckr=\R((t))$ denote the resulting real form, the field of
real formal Laurent series.
The restriction of the conjugation $c$ preserves the ring $\co=\C[ [ t] ]$ of
formal power series, and we let
$\cor=\R[ [t] ]$ denote the resulting real form, 
the ring of real formal power series.

In Section~\ref{sinv}, we extended the conjugation $\theta$
of $G$ to a conjugation $\theta$ of $G(\ck)$. 
Since the conjugation $\theta$ of $G(\ck)$ 
restricts to a conjugation of the subgroup $G(\co)$, it
induces a conjugation of the affine Grassmannian $\affgr=G(\ck)/G(\co)$.
We define the real affine Grassmannian $\raffgr$ 
to be the resulting real form of $\affgr$.
We identify it with the quotient
$$
\raffgr=G_\R(\ckr)/G_\R(\cor)
$$
where $G_\R(\ckr)$ denotes the group of $\ckr$-valued
points of $G_\R$, and
$G_\R(\cor)$ the subgroup of $\cor$-valued
points.
We shall only be interested in the structure of $\raffgr$ as a topological space, which it inherits
as a subspace of $\affgr$.
(Although we do not use the result,
it was first recongized by Quillen, then proved in~\cite{Mitch88}, that $\raffgr$ 
is homotopy equivalent to the based loop space
of the symmetric space $X_c$.)



The restrictions of the two stratifications of $\affgr$ to
the real form $\raffgr$ provide stratifications of $\raffgr$.
The strata $\rgstrat$ 
are the orbits of the action of $G_\R(\cor)$, and
the dual strata $\rgdstrat$ 
are the orbits of the action of $G_\R(\ipolyr)$.
By the Cartan decomposition, the stratum
$\rgstrat$ and dual stratum $\rgdstrat$ are non-empty
if and only if $\lambda\in\rdcowt$.
For $\lambda\in\rdcowt$, 
we define the core $\rgcore$ to be the 
intersection $\rgstrat\cap\rgdstrat$.

\begin{Proposition}\label{pmatrstrat}
The real affine Grassmannian $\raffgr$ is the disjoint union
$$
\raffgr=\bigsqcup_{\lambda\in\rdcowt} \rgstrat
\qquad
\raffgr=\bigsqcup_{\lambda\in\rdcowt} \rgdstrat.
$$
The closures of the strata and dual strata are the disjoint unions
$$
\overline\rgstrat=\bigsqcup_{\mu\leq\lambda} S^\mu_\R
\qquad
\overline\rgdstrat=\bigsqcup_{\lambda\leq \mu} T^\mu_\R.
$$
The stratum $\rgstrat$ is a smooth real quasi-projective variety,
and the core $\rgcore$ is isomorphic to the flag variety $G_\BR/P_\BR(\lambda)$.
\end{Proposition}

\begin{proof}
Thanks to Propositions~\ref{pmatstrat},~\ref{pmatdstrat}, and~\ref{pcore},
it only remains to prove the assertion~($\star$): if
$\lambda-\mu$
is a non-negative integral
linear combination of positive coroots of $G$, then
$$
S^\mu\subset\overline {S^\lambda}\mbox{ and }
T^\lambda\subset\overline {T^\mu}.
$$
Let $R^\pos\subset\cowt$ be the set of positive coroots of $G$,
and let $Q^\pos\subset\cowt$ be the semigroup generated by $R^\pos$.

\begin{Lemma}\label{lrcowtord}
The semigroup $\rcowt\cap Q^\pos$ 
is generated by elements of the form
$\alpha\in\rcowt\cap R^\pos$, and $\theta(\alpha)+\alpha\in\rcowt$,
for $\alpha\in R^\pos$.
\end{Lemma}

\begin{proof}
Let $2\ch\rho_M$ denote the sum of the positive roots of the complexification $M\subset P$ of the Levi factor $M_\BR\subset P_\BR$
of the minimal parabolic subgroup $P_\BR\subset G_\BR$.
The involution $\theta$ fixes elements of $\rcowt$ and takes $2\ch\rho_M$
 to its negative.
Therefore $\rcowt$ pairs trivially with $2\ch\rho_M$,
and so we may write $\beta\in\rcowt\cap Q^\pos$ 
uniquely as a sum $\beta=\sum_i \alpha_i$ of simple coroots $\alpha_i$
of the centralizer $L$ 
of $2\ch\rho_M$.
Since $\theta$ preserves the set of simple coroots of $L$, and the expression of $\beta$
as a sum of the simple coroots of $L$
is unique, $\theta$ also
preserves the set $\{\alpha_i\}$ of simple coroots counted with multiplicities
appearing in the sum.
The assertion follows by induction on the size of this set.
\end{proof}

By Lemma~\ref{lrcowtord}, it suffices to prove assertion ($\star$) when $\lambda-\mu$
is a positive coroot $\alpha$, or when $\lambda-\mu$ is of the form
$\theta(\alpha)+\alpha$, for a positive coroot $\alpha$
but $\lambda-\mu$ is not a multiple of a positive coroot.
In the first case, we may find $\on{SL}_2(\R)\subset G_\R$ such that $\alpha$
is its positive coroot. Let $U_\alpha\subset\on{SL}_2(\R)$ be the one-parameter subgroup
corresponding to $\alpha$. The orbit through $\lambda$ of the subgroup
$U_{\alpha}(rt)\subset G_\R(\cor)$, for $r\in\R$, 
is isomorphic to $\R$, and its closure is isomorphic to $\R\BP^1$
with $\mu$ the point at infinity. 
In the second case, we may find $\on{SL}_2(\C)\subset G_\R$ such that $\theta(\alpha)+\alpha$
is its positive coroot. 
Let $U_{\theta(\alpha)+\alpha}\subset\on{SL}_2(\BC)$ be the one-parameter subgroup
corresponding to ${\theta(\alpha)+\alpha}$.
Then the orbit through $\lambda$ of the subgroup
$U_{\theta(\alpha)+\alpha}(ct)\subset G_\R(\cor)$, 
for $c\in\C$, is isomorphic to $\C$, and its closure is isomorphic to $\C\BP^1$
with $\mu$ the point at infinity. We conclude $S^\mu\subset\overline {S^\lambda}$.
By taking the opposite subgroups with $t^{-1}$ substituted for $t$ and acting on $\mu$,
we conclude $T^\lambda\subset\overline {T^\mu}.$
\end{proof}


Recall the automorphism $\phi_z:\affgr\to\affgr$ 
induced by the automorphism $t\mapsto z^{-1}t$ of the field $\ck=\C((t))$, for $z\in\BC$.   
When $z\in\Rp$, the automorphism $t\mapsto z^{-1}t$ preserves the field $\ckr=\R((t))$, and so 
the induced automorphism of $\affgr$ preserves the real
form $\raffgr$. 
Theorem~\ref{tflow} immediately implies the following.

\begin{Theorem}\label{trmorse} 
The cores $\rgcore$ are the fixed points of the $\BR^\times$-action on $\raffgr$
given by $\phi_z$.
The strata $\rgstrat$ are the ascending spaces 
$$
\rgstrat=\{x\in\raffgr|\lim _{z\to\infty} \phi_z(x)\in \rgcore\},
$$
and the dual strata $\rgdstrat$ are the descending spaces
$$
\rgdstrat=\{x\in\raffgr|\lim _{z\to0} \phi_z(x)\in \rgcore\}.
$$
\end{Theorem}


\section{Beilinson-Drinfeld Grassmannians}\label{smatbdg}
Let $\C=\Spec(\poly)$ denote the affine line,
and $\C^\times=\spec(\lpoly)$ the punctured affine line $\C\setminus\{0\}$. 
For any $\C$-algebra $\alg$, let $\C_\alg$, respectively $\pC_\alg$, denote 
the product $\C\times\spec(\alg)$, respectively $\pC\times\spec(\alg)$.
The following proposition (first proved for $G=\on{SL}_n$ in \cite[Proposition 2.1]{BeauLa94}) 
gives an alternative characterization
of the functor $\alg\mapsto G(\alg((t)))/G(\alg[[t]])$
from $\C$-algebras to sets.

 
\begin{Proposition}[\cite{LaSo97}, Proposition 3.10]
There is a canonical isomorphism between the functor
$\alg\mapsto G(\alg((t)))/G(\alg[[t]])$ and
the functor
$$\alg\mapsto
\left(
\begin{array}{l}
  \tors \mbox{ a $G$-torsor on } \C_\alg,\\
  \nu \mbox{ a trivialization of } \tors \mbox{ over } \pC_\alg
\end{array}
\right).
$$
\end{Proposition}
The data are to be taken up to isomorphism.
We think of the trivialization as a section ${\nu:\pC_\alg\to\tors|\pC_\alg}$, so that
an isomorphism from a pair $(\tors_1,\nu_1)$ to another $(\tors_2,\nu_2)$
is an isomorphism of torsors $\phi:\tors_1\to\tors_2$ such that
$\nu_2 = \phi\circ\nu_1$.


\subsection{Definition of $\taffgr$}
For any $\C$-algebra $\alg$, let $\C(\alg)=\spec(\alg[t])$ denote the 
$\alg$-valued points of the affine line $\C$.
Consider the following 
functor from $\C$-algebras to sets
$$\alg\mapsto
\left(
\begin{array}{l}
x\in \C(\alg),\\
 \tors \mbox{ a $G$-torsor on } \C_\alg,\\
\nu \mbox{ a trivialization of } \tors \mbox{ over } 
\C_\alg \setminus (x\cup -x)
\end{array}
\right).
$$
The data are to be taken up to isomorphism.
We think of the points $x,-x:\spec(\alg)\to \C$
as subschemes of $\C_\alg$ by taking their graphs.
We think of the trivialization as a section
$\nu:\C_\alg\setminus (x\cup-x)\to\tors|{\C_\alg\setminus (x\cup-x)}$, 
so that an isomorphism from a pair 
$(\tors_1,\nu_1)$ to another $(\tors_2,\nu_2)$
is an isomorphism of torsors $\phi:\tors_1\to\tors_2$ such that
$\nu_2 =\phi\circ \nu_1$.

The functor is known to 
be represented by an ind-finite type ind-scheme. 
We refer to its space of $\BC$-points equipped with the analytic topology
as the Beilinson-Drinfeld
Grassmannian of $G$, and denote it by $\taffgr$. (See~\cite{BD,MV00} for more details.)


We have the natural projection
$\taffgr\to\C$ defined by $(x,\tors,\nu)\mapsto  x.$ 

\begin{Proposition}[\cite{BD}, 5.3.10]\label{pmatfiber}
There are canonical isomorphisms
$$
\taffgr|\{0\}\stackrel{\sim}{\to}\affgr
$$
$$
\taffgr|(\C\setminus\{0\}) \stackrel{\sim}{\to}\affgr\times\affgr\times(\C\setminus\{0\}).
$$
\end{Proposition}

Observe that the flow $\phi_z:\affgr\to\affgr$
extends to a flow $\tphi_z:\taffgr\to\taffgr$
such that $\tphi_z|\{0\}=\phi_z$ under the identification $\taffgr|\{0\}\simeq\affgr$.
The extended flow is
defined by the formula
$$
\tphi_z(x,\tors,\nu)=(z\cdot x,\tors^{-z},\nu^{-z})
$$
where $z\cdot x$ denotes multiplication of $x\in\C$ by $z\in\pC$,
$\tors^{-z}$ the pull-back of $\tors$ via multiplication by $z^{-1}$,
and $\nu^{-z}$ the pull-back of the trivialization.
We have a commutative diagram
$$
\begin{array}{ccc}
\taffgr & \stackrel{\tphi_z}{\to} & \taffgr\\
\downarrow & & \downarrow \\
\C & \stackrel{z\cdot}{\to} & \C
\end{array}
$$
where the map $z\cdot$ is multiplication by $z\in\pC$.


\subsection{The real form $\sraffgr$}

Recall that $\theta$ denotes
the conjugation of $G$ with respect to $G_\R$.
Let $c$ denote the standard conjugation of $\C$
with respect to $\R$.
We define a twisted conjugation 
of $\taffgr$
by the formula
$$
\theta_\sigma(x,\tors,\nu) = (-c(x),\tors,\nu^\theta_c)
$$
where 
the section $\nu^\theta_c$ 
is the composition
$$
\nu^\theta_c:
\C\setminus(c(x)\cup-c(x))\stackrel{c}{\to}
\C\setminus(x\cup-x)\stackrel{\nu}{\to}
\tors|\C\setminus(x\cup-x)\stackrel{f^\theta_c}{\to}
\tors|\C\setminus(x\cup-x)
$$
and the map $f^\theta_c$
is defined to be $f^\theta_c(g,z)=(\theta(g),c(z))$
with respect to any trivialization
$\tors\simeq G\times\C$. 
Although the section $\nu^\theta_c$ depends on the choice
of trivialization, the isomorphism class of the pair $(\tors,\nu^\theta_c)$
is independent of the choice.

We define the twisted real Beilinson-Drinfeld Grassmannian
$\sraffgr$ to be the real form of $\taffgr$ with respect to 
the twisted conjugation $\theta_\sigma$.

The projection $(x,\tors,\nu)\mapsto x$
takes the real form $\sraffgr$ to the imaginary axis $i\R$,
since if $x\in\C$ satisfies $-c(x)=x,$ then $x\in i\R$.

\begin{Proposition}
There are isomorphisms
$$
\sraffgr|\{0\}\stackrel{\sim}{\to}\raffgr
$$
$$
\sraffgr|(i\R\setminus\{0\})\stackrel{\sim}{\to}\affgr\times(i\R\setminus\{0\}).
$$
\end{Proposition}

\begin{proof}
The isomorphisms are the restrictions of the isomorphisms of
Proposition~\ref{pmatfiber}.
\end{proof}
 
Observe that for $z\in\Rp$, the flow $\tphi_z:\taffgr\to\taffgr$
restricts to provide a flow $\tphi_z:\sraffgr\to\sraffgr$.
We have a commutative diagram
$$
\begin{array}{ccc}
\sraffgr & \stackrel{\tphi_z}{\to} & \sraffgr\\
\downarrow & & \downarrow \\
i\R & \stackrel{z\cdot}{\to} & i\R
\end{array}
$$
where the map $z\cdot$ is multiplication by $z\in\Rp$.
Therefore the critical manifolds of the flow  
are the cores 
$\rgcore\subset\sraffgr|\{0\}$, and
the ascending spaces are the strata $\rgstrat\subset\sraffgr|\{0\}$,
for $\lambda\in\rdcowt$.

For $\lambda\in\rdcowt$,
we write $\srgdstrat\subset\sraffgr$ for
the descending space
$$
\srgdstrat=\{x\in\sraffgr|\lim _{z\to0} \tphi_z(x)\in \rgcore\}
$$
of points that contract to the core $\rgcore\subset\sraffgr|\{0\}$.
The identification $\sraffgr|\{0\}\simeq\raffgr$ restricts to an identification
$
\srgdstrat|\{0\}\simeq \rgdstrat.
$
%

\begin{Proposition}\label{prorbfiber}
The isomorphism
$
\sraffgr|\{i\}\simeq \affgr
$
restricts to an isomorphism
$
\srgdstrat|\{i\}\simeq \rorb,
$
for $\lambda\in\imgp$,
and the fiber $\srgdstrat|\{i\}$ is empty, for $\lambda\in\rdcowt\setminus\imgp$.
\end{Proposition}

\begin{proof}
We first show that all of the points of the orbit $\rorb\subset\sraffgr|\{i\}$ 
contract to 
the same core $C_\R$. We then check that the core $C_\R$ is
indeed $\rgcore$.

Let $\C\proj^1=\C\cup\{\infty\}$ be the complex projective line,
and let $\R\proj^1=\R\cup\{\infty\}$ be the real projective line.
Consider the group of rational maps $g:\C\proj^1\to G$
which have poles only at $i$ and $-i$ and take
$\R\proj^1$ into $G_\R$.
We may identify this group 
with $LG_\R\subset G(\lpoly)$. 
It acts naturally on the fiber $\sraffgr|\{i\}$. 
If a point $x\in\sraffgr|\{i\}$ contracts to a point 
$x_0\in\sraffgr|\{0\}$,
then the point $g\cdot x\in\sraffgr|\{i\}$ contracts to the point
$g(\infty)\cdot x_0\in\sraffgr|\{0\}$.
Therefore if a point $x\in\sraffgr|\{i\}$ contracts to a core $C_\R$,
then 
the point
$g\cdot x\in\sraffgr|\{i\}$ also contracts to the same core $C_\R$.

To prove the proposition, it now suffices to show that for each $\lambda\in\imgp$,
there is a point $x\in\rorb\subset\sraffgr|\{i\}$ which contracts
to the core $\rgcore\subset\sraffgr|\{0\}$.
By Lemma~\ref{lgeodinimage}, for each $\lambda\in\imgp$, there exists $x\in\blg\simeq\affgr$
such that $\ttheta^\tau(x)x=\lambda$. The $LG_\R$-orbit $\rorb\subset\affgr$ 
is by construction the orbit through $x$. 
We claim this point $x\in\rorb\subset\sraffgr|\{i\}$
contracts to $\lambda\in\rgcore\subset\sraffgr|\{0\}$.
Under the inclusion $\sraffgr\subset\taffgr$, the point 
$x\in\rorb\subset\sraffgr|\{i\}$
is represented by $(x,\theta(x))\in\affgr\times \affgr\simeq\taffgr|\{i\}$.
Consider the group of rational maps $g:\C\proj^1\to G$
which have poles only at $i$ and $-i$ and take
$\R\proj^1$ into $G_c$. We may identify this group 
with $LG_c\subset G(\lpoly)$. 
Its natural action on the fiber $\taffgr|\{i\}\simeq\affgr\times\affgr$ 
is given by the formula
$g\cdot(a,b)=(ga,\tau(g)b)$
where $\tau$ denotes the time-reversal involution of $G(\lpoly)$.
As before, if a point $x\in\taffgr|\{i\}$ contracts to a point 
$x_0\in\taffgr|\{0\}$,
then the point $g\cdot x\in\taffgr|\{i\}$ contracts to the point
$g(\infty)\cdot x_0\in\taffgr|\{0\}$.
If we take $g=x^{-1} \in LG_c$, then $g\cdot(x,\theta(x))=(1,\lambda)$.
This point clearly contracts to $\lambda\in\rgcore\subset\sraffgr|\{0\}$, so the point
$x\in\rorb\subset\sraffgr|\{i\}$ does as well.
\end{proof}


\section{Specialization functor 
$R:\catd(\affgr)\to\catd(\raffgr)$}\label{smatspec}

For an ind-scheme $Z$ of ind-finite type, 
let $\mathbf D(Z)$ be the limit category of the bounded derived categories
of finite-type subschemes of $Z$.
We use the term sheaf to refer 
to objects of $\mathbf D(Z)$
and only considered derived functors. 
If a group $H$ acts on $Z$ preserving the ind-structure,
we say a sheaf is $H$-constructible if it is constructible
with respect to a stratification by $H$-invariant subschemes.

\subsection{Definition of specialization}
Consider the diagram of ind-schemes
$$
\begin{array}{ccccccccc}
\affgr\times i\Rp  & \stackrel{}{\simeq} & 
\sraffgr|i\Rp & \stackrel{j}{\to} 
& \sraffgr|i\Rnn & \stackrel{i}{\gets} & \sraffgr|\{0\} 
& \stackrel{}{\simeq} & \raffgr \\
 & \searrow & \downarrow & & \downarrow & & \downarrow & \swarrow &\\
  && i\Rp & \to & i\Rnn & \gets & \{0\} && 
\end{array}
$$ 
Define a specialization functor $\on R:\catd(\affgr)\to\catd(\raffgr)$
by the formula
$$
\on R(\sh)= i^*{}j_*(\sh\boxtimes\C_{i\Rp} )
$$
where $\C_{i\Rp}$ denotes the constant sheaf on ${i\Rp}$.
Informally speaking, the specialization 
is the nearby cycles in the family $\sraffgr|i\Rnn\to i\Rnn$.


\subsection{Relation to $K(\ck)$-orbits}

\begin{Theorem}\label{tconst}
The specialization $\on R:\catd(\affgr)\to\catd(\raffgr)$ 
takes $K(\co)$-constructible
objects to $G_\R(\cor)$-constructible objects.
\end{Theorem}

\begin{proof}
The intersection of the subgroups $K(\co)$ and $G_\R(\cor)$ in the group
$G(\co)$ is the subgroup $K_c(\cor)$ of $\cor$-valued points of 
the compact group $K_c$.

\begin{Lemma}
The orbits of $K_c(\cor)$ in 
$\raffgr$ coincide with the orbits of $G_\R(\cor)$.
\end{Lemma}

\begin{proof}
First note that the inclusion $K_c\subset G_\BR$ induces an isomorphism of component groups.
Since the automorphism $t\mapsto rt$,
for $r\in\Rp$, induces a contraction of $K_c(\cor)$ to $K_c$ and of $G_\BR(\cor)$ to $G_\BR$,
the inclusion $K_c(\cor)\subset G_\BR(\cor)$ also induces an isomorphism of component groups.
Moreover, thanks to this contraction,
it now suffices to prove an infinitesimal version
of the result at a core point $x\in\rgcore\subset\raffgr$.
In other words, it suffices to show that the Lie algebra $\liek_c(\cor)$
of the group $K_c(\cor)$ surjects onto the tangent space $T_x\rgstrat$
where the stratum $\rgstrat$ is the $G_\R(\cor)$-orbit through $x$.
Equivalently, it suffices to show
the Lie algebra $\liek_c(\cor)$ is transverse
to the kernel $\EuFrak s_x$ of the surjection $\liegr(\cor)\to T_x\rgstrat$.
For each $k$, we may identify the subspace $t^k\liegr\subset\liegr(\cor)$
with the Lie algebra $\liegr$. 
The intersection $\liek_c(\cor)\cap t^k\liegr$ is a compact
subalgebra of $\liegr$, and $\EuFrak s_x\cap t^k\liegr$ contains a parabolic
subalgebra of $\liegr$. 
Thus in each subspace $t^k\liegr\subset\liegr(\cor)$ 
the subalgebras $\liek_c(\cor)$ and $\EuFrak s_x$
are transverse, and the result follows.
\end{proof}

By the lemma, it suffices to show that the specialization
takes $K(\co)$-constructible
objects to $K_c(\cor)$-constructible objects.

Consider the group $\lockc$ of maps $\spec(\R[t]_{(t)})\to K_c$
where $\R[t]_{(t)}$ denotes the local ring of rational functions on $\R$
without poles at $\{0\}\in\R$. The group $\lockc$ naturally acts on the
germ of $\sraffgr|\{0\}$ inside $\sraffgr$. 
Therefore the following lemma implies the theorem.

\begin{Lemma}
The orbits of $\lockc$ in 
$\raffgr\simeq\sraffgr|\{0\}$ coincide with the orbits of $K_c(\cor)$.
\end{Lemma}

\noindent{\it Proof.}
For any orbit of $K_c(\cor)$, there exists $k$
such that the action of $K_c(\cor)$ on the orbit descends
to an action of the group $K_c(\jet^k_\R)$ where $\jet^k_\R=\cor/(t^k)\cor$.
Now it is not difficult
to prove that the natural map $\lockc\to K_c(\cor)\to K_c(\jet^k_\R)$
is surjective. See~\cite[Lemma 4]{Gait01} for details in a similar
situation.
\end{proof}
%


\subsection{Relation to $LG_\R$-orbits}
Recall that for the flow $\tphi_z:\sraffgr\to\sraffgr$, 
and for $\lambda\in\rdcowt$,
we have the critical manifolds $C^\lambda_\R$,
the ascending spaces $S^\lambda_\R$, 
and the descending spaces $\tilde T^\lambda_\R$,
For $\lambda\in\rdcowt$, we write 
$$
s_\lambda:S^\lambda_\R\to\sraffgr\qquad
\tilde t_\lambda:\tilde T^\lambda_\R\to\sraffgr
$$
for the inclusions, and
$$
c_\lambda:S^\lambda_\R\to C^\lambda_\R\qquad
\tilde d_\lambda:\tilde T^\lambda_\R\to C^\lambda_\R
$$
for the contractions.

%

\begin{Theorem}\label{tcalc}
For any object $\sh\in\catd(\sraffgr)$ that is 
$\Rp$-constructible with respect to the flow
$\tphi_z:\sraffgr\to\sraffgr$, for $z\in\Rp$, we have an isomorphism
$$
(c_\lambda)_! (s_\lambda)^*\sh\simeq (\tilde d_\lambda)_* (\tilde t_\lambda)^!\sh,
\mbox{ for any } \lambda\in\rdcowt.
$$
\end{Theorem}
%

\begin{proof}
The proof is an application of \cite[Proposition 9.2]{GM93}
whose terminology we adopt.
We may restrict our attention to the support of $\sh$
which is a finite-dimensional variety.
Then by \cite{Nad02}, we may work in an ambient smooth
variety to which the automorphism $\phi_z$, for $z\in\BR^\times$, extends.
Furthermore, any point of the ambient smooth variety 
has a neighborhood homeomorphic
to an open set in $\BR^n$, 
such that under this homeomorphism $\phi_z$ acts linearly on $\BR^n$ via a set of
integral weights $w_1,\ldots,w_n$.
We may assume that $w_k > 0$, for $k\leq l$, $w_k=0$, for $l<k<m$, and $w_k < 0$, for $m\leq k$. 
Define an indicator map $\R^n\to \R^2$ by
 $$
(x_1,...,x_n) \mapsto 
(x_1^{N/w_1} + ... + x_m^{N/w_l}, 
                        x_l^{N/w_m} + ... + x_n^{N/w_n})
$$
where $N$ is the least common multiple of $w_1,\ldots, w_n$.
This is an equivariant map where $\BR^\times$ acts on
$\BR^2$ by $z\cdot (x,y) = (z^N x, z^{-N} y)$.  Therefore there is a stratification of $\BR^2$
such that the indicator map is a stratified 
map, and there are no non-open strata in a neighborhood of
$(0,0)$
other than $(0,0)$ itself.
Now \cite[Proposition 9.2]{GM93} applies and gives that the two sheaves
of the theorem have isomorphic hypercohomology.
In fact, the same arguments apply to the hypercohomology of any restrictions of the sheaves.
Thus by induction on strata,
we conclude that the sheaves themselves are isomorphic.
%
\end{proof}

Recall that by Proposition~\ref{prorbfiber}, we have an isomorphism
$\tilde T^\lambda_\R|\{i\}\simeq\rorb$, for $\lambda\in\imgp$,
and the fiber $T^\lambda_\R|\{i\}$ is empty, for $\lambda\in\rdcowt\setminus\imgp$.
For $\lambda\in\imgp$,
we write 
$$
t_\lambda:\rorb\to\sraffgr
$$
for the restriction of the inclusion $\tilde t_\lambda:\tilde T^\lambda_\R\to\sraffgr$,
and 
$$
d_\lambda:\rorb\to C^\lambda_\R
$$
for the restriction of the contraction
$\tilde d_\lambda:\tilde T^\lambda_\R\to C^\lambda_\R$.

\begin{Corollary}\label{chyp}
For any object $\sh\in\catd(\affgr)$ that is 
$\Rp$-constructible with respect to the flow
$\phi_z:\affgr\to\affgr$, for $z\in\Rp$, we have 
$$
(c_\lambda)_!(s_\lambda)^*\on R(\sh)\simeq (d_\lambda)_* (t_\lambda)^!\sh
\quad\mbox{ if } \lambda\in\imgp 
$$
$$
(c_\lambda)_!(s_\lambda)^*\on R(\sh)=0 \quad\mbox{ if } \lambda\in\rdcowt\setminus\imgp. 
$$
\end{Corollary}

\begin{proof}
Note that the constructibility of $\sh$, 
implies the constructibility
of $\sh\boxtimes\C_{i\Rp}$ required in the hypothesis of Theorem~\ref{tcalc}.
For the first assertion, we calculate
$$
\begin{array}{ccll}
(c_\lambda)_!(s_\lambda)^*\on R(\sh) & = & 
(c_\lambda)_!(s_\lambda)^*i^*j_*(\sh\boxtimes\C_{i\Rp})
& \mbox{(definition of $\on R$)} \\
& \simeq &  (c_\lambda)_!(s_\lambda)^* j_*(\sh\boxtimes\C_{i\Rp})
& \mbox{(composition of maps)} \\
& \simeq &   (\tilde d_\lambda)_* (\tilde t_\lambda)^!j_*(\sh\boxtimes\C_{i\Rp})
& \mbox{(Theorem~\ref{tcalc})} \\
& \simeq &   (\tilde d_\lambda)_* j_*(\tilde t_\lambda)^!(\sh\boxtimes\C_{i\Rp})
& \mbox{(base change)} \\
& \simeq &   (\tilde d_\lambda)_*(\tilde t_\lambda)^!(\sh\boxtimes\C_{i\Rp})
& \mbox{(composition of maps)} \\
& \simeq &   (d_\lambda)_*(t_\lambda)^!\sh
& \mbox{(integration along $\Rp$-orbits)}. \\
\end{array}
$$
The second assertion follows from the first if 
for $\lambda\in\rdcowt\setminus\imgp$,
we set $\rorb$ to be the empty set.
\end{proof}


\section{Poset structure: closure relations}\label{scr}

In this final section, we confirm that the poset structures on the sets
of orbits of the groups $K(\ck)$ and $LG_\R$ in the affine Grassmannian
$\affgr$ are as asserted in the introduction. 

Recall that $\imgp\subset\rdcowt$ is the semigroup
of dominant real coweights mapped by $[\cdot]:\rdcowt \to \pi_1(X_c)$ 
into the image of $\pi_*:\pi_1(G)\to\pi_1(X)$. 
We order $\imgp\subset \rdcowt$ by restricting the natural order from 
the dominant coweights:
$\lambda\leq \mu$ if and only if
$\mu-\lambda$ is a non-negative integral linear combination of positive
coroots of $G$.

\begin{Theorem}
The bijection $K(\ck)\backslash \affgr\risom\imgp$
is a poset isomorphism, 
and the bijection $LG_\R\backslash \affgr\risom\imgp$ 
is a poset anti-isomorphism.
\end{Theorem}

\noindent{\it Proof.}
Recall that we call two dominant real coweights $\lambda\leq \mu\in\rdcowt$
primitively related 
if there is no dominant
real coweight $\nu\in\rdcowt$ with $\lambda<\nu<\mu$.
By Proposition~\ref{pposetbounds}, Proposition~\ref{ppriminter},
and the fact that the relation of orbit closure is transitive, 
it suffices to show 
the following.

\begin{Proposition}
If $\lambda\leq \mu\in\imgp$ are primitively related,
then the intersection $\orb_K^\mu\cap\rorb$ is nonempty.
\end{Proposition}

\noindent{\it Proof.}
By Lemma~\ref{lgeodinimage}, there is a point $x\in\orb_K^{\mu}\cap\orb_\R^{\mu}\subset\affgr$.
Consider the $K(\co)$-orbit $S\subset\affgr$ through $x$.
We shall show that $S\cap\rorb$ is non-empty
by applying the results of the previous section.
Let $\C_{S}$ be the $*$-extension to $\affgr$ of
the constant sheaf on $S$.
By Corollary~\ref{chyp}, we have an isomorphism
$$
(c_\lambda)_!(s_\lambda)^*\on R(\C_S)\simeq (d_\lambda)_* (t_\lambda)^!\C_S.
$$
If the intersection $S\cap\rorb$ were empty, then the sheaf $(t_\lambda)^!\C_S$
would vanish.
Therefore the following lemma implies the proposition and thus the theorem as well.

\begin{Lemma}
The sheaf $(c_\lambda)_!(s_\lambda)^*\on R(\C_S)$ does not vanish.
\end{Lemma}

\begin{proof}
By Theorem~\ref{tconst}, the sheaf $\on R(\C_S)$ is $G_\R(\cor)$-constructible,
thus it suffices to show $(s_\lambda)^*\on R(\C_S)$ does not vanish since
$c_\lambda$ has affine fibers.
Simplifying the expression
$$
\begin{array}{ccll}
(s_\lambda)^*\on R(\C_S) & = & 
(s_\lambda)^*i^*j_*(\C_S\boxtimes\C_{i\Rp})
& \mbox{(definition of $\on R$)} \\
& \simeq & (s_\lambda)^* j_*(\C_S\boxtimes\C_{i\Rp})
& \mbox{(composition of maps)}, \\
\end{array}
$$
we see that it suffices to show $S^\lambda_\R\subset\sraffgr|\{0\}$
is in the closure of $S\times i\Rp\subset\sraffgr|i\Rp$.

Let $\overline{S_{i\Rp}}\subset\sraffgr$ be the closure
of
$S\times i\Rp \subset\sraffgr|i\Rp$. 
By similar arguments to those given in the proof of Proposition~\ref{prorbfiber}, 
$\overline{S_{i\Rp}}$ 
contains the core $C^{\mu}_\R\subset\sraffgr|\{0\}$ since the intersection 
$S\cap\orb_\R^{\mu}\subset\sraffgr|\{i\}$ is nonempty. 
Thus by Proposition~\ref{pmatrstrat} and Theorem~\ref{tconst}, 
$\overline{S_{i\Rp}}$ 
contains the strata $S^\nu_\R\subset\sraffgr|\{0\}$, for all $\nu\leq\mu$. In particular,
it contains the stratum $S^\lambda_\R\subset\sraffgr|\{0\}$. 
\end{proof}


\bibliographystyle{alpha}
\bibliography{ref}


\end{document}